\documentstyle[12pt,psfig,amssymb]{amsart}

\newcommand{\new}{\newcommand}
\new{\tnsr}{\otimes}
\new{\tensor}{\otimes}
\new{\superset}{\supset}
\new{\union}{\cup}
\new{\CC}{{\Bbb C}}
\new{\NN}{{\Bbb N}}
\new{\RR}{{\Bbb R}}
\new{\ZZ}{\mbox{$\Bbb Z$}}
\new{\FF}{{\Bbb F}}
\new{\TT}{{\Bbb T}}
\new{\tr}{\operatorname{tr}}
\new{\abs}[1]{\left|#1 \right|}
\new{\norm}[1]{\left\| #1 \right\|}
\new{\bracket}[1]{\langle #1 \rangle}
\new{\defequals}{\stackrel{\rm def}{=}}
\new{\comb}[2]{{#1 \choose #2}}
\new{\lbl}[1]{\label{#1}
        \if\draft y
                \smash{\makebox[0pt]{\hspace{-0.5in}
                        \raisebox{8pt}{\rm\tiny #1}}}
        \fi
}

\newtheorem{proposition}{Proposition}
\newtheorem{theorem}{Theorem}

\theoremstyle{definition}
\newtheorem{lemma}{Lemma}
\newtheorem{claim}{Claim}
\newtheorem{corollary}{Corollary}
\newtheorem{conjecture}{Conjecture}
\newtheorem{definition}{Definition}

\theoremstyle{remark}
\newtheorem{remark}{Remark}

\newcounter{letter}
\newcounter{numeral}
\newcounter{Numeral}
\newenvironment{alphalist}{
\begin{list}{{\normalshape(\alph{letter})}}{\usecounter{letter}}
}{\end{list}}

\newenvironment{Romanlist}{
\begin{list}{\Roman{Numeral}}{\usecounter{Numeral}}
}{\end{list}}

\new{\prop}[1]{\begin{proposition} \lbl{#1}}
\new{\thm}[1]{\begin{theorem} \lbl{#1}}
\new{\lem}[1]{\begin{lemma} \lbl{#1}}
\new{\clm}[1]{\begin{claim} \lbl{#1}}
\new{\cor}[1]{\begin{corollary} \lbl{#1}}
\new{\cnj}[1]{\begin{conjecture} \lbl{#1}}
\new{\fig}{\begin{figure}[hbt]}
\new{\eq}[1]{\begin{equation} \lbl{#1}}

\new{\eprop}{\end{proposition}}
\new{\ethm}{\end{theorem}}
\new{\elem}{\end{lemma}}
\new{\eclm}{\end{claim}}
\new{\ecor}{\end{corollary}}
\new{\ecnj}{\end{conjecture}}
\new{\efig}[2]{\caption{#1} \label{#2} \end{figure}}
\new{\eeq}{\end{equation}}

\new{\pic}[5]{\raisebox{#3pt}{
\hspace{#4pt}\psfig{file=#1.ps,height=#2pt,silent=}\hspace{#5pt}}}
\new{\chr}[1]{\mathchoice{
\pic{#1}{12}{-3}{-1}{2}}{
\pic{#1}{12}{-2}{-3}{2}}{
\pic{#1}{9}{-2}{-3}{1}}{
\pic{#1}{7}{-1}{-1}{0}}}

\begin{document}
\author{Stephen Sawin}
\address{Math Dept. 2-265, MIT,
Cambridge, MA 02139-4307, (617) 253-4995}
\email{sawin@@math.mit.edu}

\new{\qtr}{\operatorname{qtr}}
\new{\qdim}{\operatorname{qdim}}
\new{\dom}{\operatorname{dom}}
\new{\codom}{\operatorname{codom}}
\new{\comp}{\circ}
\new{\End}{\operatorname{End}}

\new{\draft}{n}

\title[TQFT's]{Links, Quantum Groups and TQFT's}
\thanks{This research supported in part by NSF postdoctoral Fellowship
       \#23068 }

\keywords{cobordisms, Hopf algebras, Jones polynomial, knots, link
  invariants,
  quantum groups, tangles, three-manifolds, topological quantum field
  theory}
\subjclass{57M25; Secondary 16W30, 57M30}

\begin{abstract}
The Jones polynomial and the Kauffman bracket are constructed, and
their relation with knot and link theory is described.  The quantum
groups and tangle functor formalisms for understanding these
invariants and their descendents are given. The quantum group
$U_q(sl_2)$, which gives rise to the Jones polynomial, is constructed
explicitly.  The $3$-manifold invariants and the axiomatic topological
quantum field theories which arise from these link invariants at
certain values of the parameter are constructed.
\end{abstract}
\maketitle

\section*{Introduction}

In studying a class of mathematical objects, such as knots, one
usually begins by developing ideas and machinery to understand
specific features of them.   When invariants are discovered, they
arise out of such understanding and generally give
information about those features.  Thus they are often immediately
useful for proving theorems, even if they are difficult to compute.

This is not the situation we find ourselves in with the invariants of
knots and $3$-manifolds which have appeared since the Jones
polynomial.  We have a wealth of invariants, all readily computable,
but standing  decidedly outside the traditions of knot and
$3$-manifold theory.  They lack a geometric interpretation, and
consequently have been of almost no use in answering questions one
might have asked before their creation.  In effect, we are left
looking for the branch of mathematics from which these should have
come organically.

In fact, we have an answer of sorts.  In \cite{Witten89a},  Witten
gave a heuristic definition of the Jones polynomial in terms of a topological
quantum field theory, following the outline of a program proposed by
Atiyah \cite{Atiyah88}.  Specifically, he considered a knot in a
$3$-manifold  and a
connection $A$ on some principal $G$-bundle, with $G$ a simple Lie
group.  The Chern-Simons functional associates a number $CS(A)$ to
$A$, but it is well-defined only up to an integer, so the quantity
$\exp(2\pi i k CS(A))$ is well-defined.
Also, the holonomy of the connection around the knot is an element of
$G$ well-defined up to conjugation,
so the trace of it with respect to a given representation is
well-defined.   Multiplying these two gives a number depending on the knot,
the manifold, the representation and the
connection.  The magic
comes when we average over all connections and all principal
bundles:  Of course, this makes no
sense, since there is no apparent measure on the infinite-dimensional
space of connections.  But proceeding heuristically, such an average should
depend
only on the manifold, representation and the isotopy type of
the knot.  Witten  argued using a close correspondence with conformal
field theory, that when the manifold is $S^3$
and the representation the fundamental one, this
invariant had
combinatorial properties that forced it to be the analogue of the
Jones polynomial for the given group
(Specifically, the Jones polynomial at certain values of the
parameter).
Needless to say, a long physics tradition of very successful heuristic
reasoning along these lines suggested to Witten that this ill-defined
average should make sense in this case.

The fundamental problem of the field, then, is to develop this
nonrigorous but beautiful geometric interpretation of the new
invariants into a rigorous one that reasonably captures its flavor.
This is important for at least two reasons.

First, such an interpretation of the invariants is probably necessary
for their application to topological questions.  While the invariants
have yet to prove
themselves in this regard, it is to be hoped that as part of a
well-developed theory they could go a long way towards unlocking the
mysteries of low-dimensional topology.  This theory seems most likely
to arise out of something like Witten's construction, which relates the
invariants explicitly to  gauge theory and other geometry.

Second, Chern-Simons field theory is
 a field theory which is clearly
nontrivial but which admits an exact combinatorial solution.  In this
sense, interpreting Witten's work rigorously is just part of a large
endeavor within mathematical physics to understand mathematically the
whole of quantum field theory.  Such an understanding seems far more
tractable for the Chern-Simons theory, with its combinatorial,
finite-dimensional expression, than for genuinely physical theories.
{}From this point of view topology may be seen as a laboratory for a
particularly simple kind of physics, the understanding of which may
point us towards a rigorous foundation for the more complex and
physically interesting theories.

Work towards this goal proceeds along two fronts.  On the one hand are
efforts to
understand the geometry and physics of the Chern-Simons field theory,
most notably recent efforts to do perturbation theory in this context
\cite{AS92,AS94,Kontsevich??}.  On the other are combinatorial,
algebraic, and topological efforts to understand the invariants
themselves, especially with an eye towards structures suggested by the
physics.  This is the subject of the present article.

We begin in Section 1 with a discussion of the Jones polynomial, or
more precisely the Kauffman bracket form of it.  We mention other
invariants, but throughout our focus will be on this example as the
easiest.  Section 2  first gives a functorial framework for viewing
these invariants, constructing the functor explicitly for the Kauffman
bracket, and then shows
such a functorial setup arises naturally from a certain algebraic
structure: a ribbon Hopf algebra.  Section 3 discusses quantum groups,
the family of
ribbon Hopf algebras associated to each Lie algebra.  We construct
the quantum group associated to $sl_2$ explicitly, and show that with the
fundamental representation it gives the Kauffman bracket as its
invariant.

The last four sections give the construction of the topological quantum
field theories.  These bear further comment.  Witten's construction
obviously gives a $3$-manifold invariant as well, by considering an
arbitrary manifold with the empty link.  This formed some of the inspiration
for Reshetikhin and Turaev's construction of the $3$-manifold
invariants \cite{RT91} using the algebraic machinery of quantum
groups.  But Atiyah observed that this averaging
process (called path integration by physicists), if taken at its word,
implied some strong statements about how the invariants behave under
cutting and pasting.  These statements were formulated by Atiyah into
his axioms for a topological quantum field theory, or TQFT
\cite{Atiyah89,Atiyah90b}.  It was shown in
\cite{Walker??} and \cite{Turaev94} that the $3$-manifold invariants
satisfy these axioms, although the key ideas of this proof already
appeared in \cite{RT91}, and in a less rigorous form in
\cite{Witten89b,Witten89a}.  One of the principal aims of this paper
is to
give a simple proof of this fact, which is the best rigorous
connection we have between these invariants and actual quantum field
theory.

Section 4 discusses the properties of the quantum groups at roots of
unity, summarized as their being {\em modular Hopf algebras}, which
allow one to construct TQFT's from them.  Section 5 gives a categorical
formulation of TQFT's.  Section 6 sketches a purely combinatorial
description of the relevant category, that of biframed $3$-dimensional
cobordisms, in terms of surgery on links.  This allows us in Section
7 to construct the TQFT out of the link invariants we have already
constructed.

This paper is intended to introduce and invite a large
mathematical audience to this field.  The focus is on giving the
flavor and illustrating some of the power of these simple ideas.  This
is attempted by proving an important result with a minimum of
machinery, rather than by surveying the whole of the field, or giving
a full account of the  machinery.  Consequently many interesting
areas go unmentioned or barely mentioned.
By the
same token, the references emphasize pointing the interested reader to
good introductory accounts, rather than a thorough assignment of
credit.  I apologize for any omissions, and point the reader to those
same introductory works, many of which have excellent
bibliographies, for the complete story.

I would like to thank Scott Axelrod, John Baez, Joan Birman, Dana
Fine, Vaughan
Jones, Robert Kotiuga, Haynes Miller, Andrew Pressley, Justin
Roberts, Lisa Sawin, Isadore Singer, Jim Stasheff, Vladimir Turaev,
David Vogan and Eric
Weinstein for helpful conversations, comments and suggestions.
\section{Links and Their Invariants}

The study of knots and links  begins with some
physics which, while a little
eccentric sounding to modern ears,  and
decidedly wrong, bears a striking resemblance to the physics where this
story ends.  In 1867 Lord Kelvin proposed that atoms were knotted
vortices of ether, and molecules were linked atoms
\cite{Thompson67}.  Efforts began to
solve the basic problem of knot theory:  When are two knots
the same?  Early efforts to approach this by careful tabulating and
naive searching for invariants failed completely.  Modern topology has
proved more successful.  The  approach of focusing exclusively on
the knot complement and almost exclusively on its fundamental group
has, by dint of great effort, produced a complete but thoroughly
impractical algorithm for determining if two prime knots (knots which
cannot be cut into two smaller knots) are the same
up to orientation of the knot and space \cite{Hemion92}.  The story
for links is more complicated, and not as well understood.

The Jones polynomial, with its descendants, grew from an entirely
different branch of mathematics.
Its origins bear telling as a study in serendipity,
although we will only sketch  Jones' construction, which is buried in
the construction in Section 2.  Jones, in the course of
proving an important result
about the ways in which certain algebras of operators  sit inside
each other \cite{Jones83}, constructed an algebra with a trace on it.  He
then noticed  that this algebra gave a representation of
the braid group (the mapping-class group of the plane with finitely
many points removed).  Now by Markov's
theorem \cite{BZ85}, links can be
represented by elements of the braid group, with two braids giving the
same link exactly when they can be connected by conjugation and
another move, called stabilization.   Any trace is invariant under
conjugation, and this particular trace could easily be normalized so as
to be invariant under stabilization.  Thus presenting the link as
a braid, representing the braid in the algebra, and taking the trace
gives an invariant of the link.  Since the algebra and the trace
depend on a parameter $t$, so does the invariant, and it turns out
the invariant is essentially a Laurent polynomial in $t$: the
celebrated Jones polynomial.

This invariant almost certainly does not come from  classical knot theory:  It
distinguishes links with diffeomorphic complement and detects mirror
images for example.  It is eminently
computable:  Although the algorithm is roughly exponential in the number of
crossings, and the problem is known to be $\#$P-hard \cite{JVW90},
a clever high school student can  compute it
easily.  The collection
of these invariants is a fairly
good distinguisher of knots, though it does not distinguish all knots,
and it is not even known whether it can distinguish a
knot from the unknot \cite{Garoufalidas??}.  Maddeningly, these
invariants have no geometric
interpretation (at least, no rigorous one), and consequently they have
provided few purely topological results (chiefly the
Tait conjectures, about so-called alternating knots).  With such
mysteries before us, it is time to do some  mathematics.

A link is a smooth embedding of several copies of $S^1$ into oriented
$S^3$.  A knot is a link with one component.  Links are equivalent if
there is an
orientation-preserving diffeomorphism of $S^3$ taking one to the other.  This
notion is
equivalent to the corresponding PL notion, and to the corresponding
topological notion if we restrict to  `tame' knots \cite{BZ85}.  An oriented
link
has an orientation on each component, and a framed link comes equipped
with a nonzero section of the normal bundle, also up to positive
diffeomorphism.  We will draw links as in Figure \ref{fg:linkexamples},
with only transverse-double-point crossings.  Oriented links will be
drawn with arrows indicating orientation, and framed links will be drawn
so that the framing lies in the plane  (thus in Figure
\ref{fg:linkexamples}, B represents a different framed knot from A,
because the
twist cannot be undone).
The unknot is the knot
which bounds a disk, and has an obvious preferred framing (shown in
A).  This is all
the knot theory we will need to define the Jones polynomial, using a
construction due to Kauffman \cite{Kauffman87}.

\fig
{\small
\begin{tabular*}{\textwidth}{c@{\extracolsep{\fill}} c@{\extracolsep{\fill}}
    c@{\extracolsep{\fill}} c@{\extracolsep{\fill}} c@{\extracolsep{\fill}}
    c@{\extracolsep{\fill}}c}

A. Unknot & B. Framed & C. Rt. & D. Left & E. Rt. or. & F. Left or. &
G. Hopf\\
& unknot & trefoil & trefoil & trefoil & trefoil & link\\
\pic{Link-u/unknot}{25}{-15}{0}{0}&
\pic{Link-u/pframedunknot}{30}{-15}{0}{0}&
\pic{Link-u/trefoil}{30}{-15}{0}{0}&
\pic{Link-u/othertrefoil}{30}{-15}{0}{0}&
\pic{Link-o/trefoil}{30}{-15}{0}{0}&
\pic{Link-o/othertrefoil}{30}{-15}{0}{0}&
\pic{Link-o/hopflink}{25}{-15}{0}{0}
\end{tabular*}
}
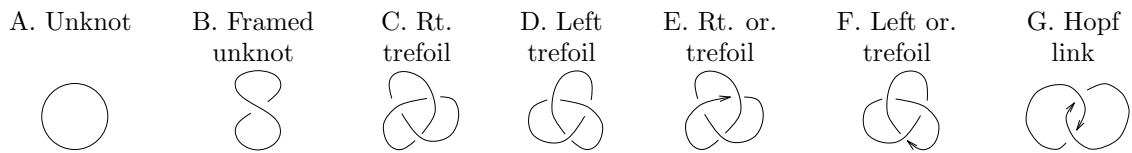
\efig{Some knots and links}{fg:linkexamples}

Let ${\cal L}$ be a projection of a framed, unoriented link $L$.
The Kauffman bracket of ${\cal L}$, $\bracket{{\cal L}}$, is an
 element of $\ZZ[A,A^{-1}]$, with $A$  an indeterminate, computed
 by the following skein relations:
\begin{align}
\bracket{\phi}&=1 \label{eq:Kempty}\\
\bracket{\chr{Dfrag-u/overcrossing}}&= A\bracket{\chr{Dfrag-u/lzero}}
+ A^{-1}\bracket{\chr{Dfrag-u/linfinity}}
\label{eq:Kskein}\\
\bracket{\chr{Dfrag-u/unknot}}&= (-A^2 -
A^{-2})\bracket{\chr{Dfrag-u/empty}}\label{eq:Kunknot}
\end{align}
where $\phi$ is the link with no components.  Equation
(\ref{eq:Kskein}), for example, says that any
time you can find three different link projections   which look
exactly the same except in a small disk,
where they look as shown in the equation, then their brackets satisfy
this equation.  Of course, this means if you happen to know the
brackets of the two projections on the right side, this tells you
the bracket of the left side.  Equation (\ref{eq:Kunknot}) is
interpreted similarly, and gives the effect of removing an unlinked
unknot from a link.
\thm{th:Kdefinition} \cite{Kauffman87}
The bracket of every projection is uniquely determined by these three
conditions, and one can compute it explicitly from them.
\ethm

\begin{pf}
By induction on the number of crossings.
\end{pf}

For example,
\begin{align*}
\bracket{\chr{Link-u/plusframedunknot}}&= A\bracket{\chr{Link-u/twounknots}}+
A^{-1}\bracket{\chr{Link-u/flatunknot}}\\
& = A(-A^2-A^{-2})^2 + A^{-1}(-A^2-A^{-2})=
 -A^{3}(-A^2-A^{-2}).
\end{align*}
More generally,
$$
\bracket{\chr{Dfrag-u/postwist}}= A
\bracket{\chr{Dfrag-u/cupandunknot}} +
A^{-1}\bracket{\chr{Dfrag-u/vase}}
=-A^{3}\bracket{\chr{Dfrag-u/cup}},
$$
i.e., adding a twist  multiplies any bracket by
$-A^{3}$.  This is actually a positive twist, corresponding
to a clockwise rotation of the framing.  A negative twist would
multiply the bracket by $-A^{-3}$. The reader can compute
that the Kauffman bracket of the right-handed trefoil (C in  Figure
\ref{fg:linkexamples}) is $A^7 + A^3 + A^{-1} - A^{-9} $.

The point is that the Kauffman bracket does not depend on
the projection, but only on the unoriented framed link.  This
remarkable fact  follows from  simple calculations, and the
framed Reidemeister's  theorem,  which says that two
projections correspond to the same framed unoriented links if and only if they
can be connected by a sequence of a certain finite set of local moves
which are roughly Moves I, II, III of Theorem \ref{th:presentation} in Section
2 (see \cite{BZ85} for
Reidemeister's theorem, the framed version follows from
\cite{Trace83}).

\thm{th:Kinvariance} \cite{Kauffman87} The Kauffman bracket takes the same
  value on two projections of the same framed link $L$.
  Henceforth refer to this value as $\bracket{L}$. \qed
\ethm

We would like a link invariant, rather than a framed link invariant,
if for no other reason than that that is what has been studied since
Kelvin (although a case can be made that framed links are a more
`natural' object to study: see e.g., Kirby's theorem in Section 5).
We have already seen that adding a positive or negative twist to a
link changes the bracket by $-A^{\pm 3}$.  It also changes the framing
by one full turn.  The idea is to correct the bracket so that this move
leaves it unchanged, without losing the framed link invariance.   To
do this, consider an {\em oriented\/} link projection, and label each crossing
as positive or negative according to whether it looks like the first or
second crossing appearing in Equation (\ref{eq:Vskein}):
I.e., according to whether one strand rotates clockwise or
counterclockwise around the other.  Then the {\em writhe} of an
oriented link
projection is  the number of positive crossings minus the
number of negative crossings.  For example, the writhe of E
 in Figure \ref{fg:linkexamples} is $+3$, that of F, its
mirror image, is $-3$.  This  is an invariant of
oriented framed links, which increases by one when a full twist is
added.  Thus if ${\cal L}$ is a projection of an oriented link $L$, the
quantity
\begin{equation}
(-A)^{-3\operatorname{writhe}({\cal L})}\bracket{{\cal L}}
\label{eq:ABJones}, \end{equation}
where $\bracket{{\cal L}}$ is understood to be the bracket of ${\cal L}$
with the orientation removed, depends only on the link $L$.

Even better, notice that in the bracket of the
trefoil, exponents increased in steps of four.  This is true in
general.  In fact, a pretty induction on the number of crossing shows
that the quantity  (\ref{eq:ABJones}) is an element of
$\ZZ[A^4,A^{-4}]$ if there are an even number of components and is
$A^2$ times such an element if the number of components is odd.  This
motivates
\begin{definition}  The quantity (\ref{eq:ABJones}), with $t^{-1/4}$
  substituted for $A$, is called the Jones polynomial of $L$, or
  $V_L(t)$.  It sends oriented, unframed links to polynomials in $t$ and
  $t^{-1}$ (times  $t^{1/2}$ if there are an odd number of components)
  and satisfies the skein relations
\begin{align}
V_\emptyset(t) &= 1 \label{eq:Vempty}\\
t^{-1}V_{\chr{Dfrag-o/oplus}}(t) - t V_{\chr{Dfrag-o/ominus}}(t)&=
(t^{1/2}-t^{-1/2})V_{\chr{Dfrag-o/ozero}}(t)
\label{eq:Vskein}
\end{align}
which uniquely determine it.
\end{definition}

One can use Equations (\ref{eq:Vempty}) and (\ref{eq:Vskein}) to
compute the Jones polynomial more efficiently than via the Kauffman
bracket, but it is a bit trickier (see, e.g., \cite{LM88}).  We should
note that the Jones polynomial as it often appears is our Jones
polynomial divided by $-t^{1/2}-t^{-1/2}$, so that its value on the
unknot is $1$.

It is worth doing a small computation to see what we have.  From
the calculations above,  the Jones
polynomial of the right-handed trefoil, E in Figure
\ref{fg:linkexamples}, is $t^{1/2}(t^4-t^2-t-1)  $.  But from
the definition, the Jones
polynomial of the mirror image of a link (a projection of which is
gotten by reversing every crossing of a projection of the link) has
the same Jones polynomial, except with $t$ replaced by $t^{-1}$.  Thus
the left-handed trefoil, F in Figure \ref{fg:linkexamples},
has invariant $t^{-1/2}(t^{-4}-t^{-2}-t^{-1} -1)  $.  Since these are
different, the two trefoils are
not equivalent!  Already we can see there is more afoot here than
simply the fundamental group of the complement, the major tool of
classical knot theory,  which must be the same for both.  In
fact, it was not until 1914 that Dehn proved the inequivalence of
these two.

The Jones polynomial was quickly generalized by a host of people
\cite{HOMFLY85,PT87} to a similar,  two-variable polynomial, now known as the
HOMFLY polynomial (the name, formed from their initials,  was
unfortunately coined before the work
of Przyticki and Traczyk was well known), and by Kauffman  to
another two-variable polynomial,  the Kauffman polynomial,
with very similar properties \cite{Kauffman85,Kauffman87}.

  Connections with  statistical mechanics were quickly noticed
  \cite{Jones85,Kauffman88}.  Many exactly-solvable statistical
  mechanical systems had recently been constructed \cite{Baxter82}
  from solutions of
  the quantum Yang-Baxter Equation (\ref{eq:YB}).  Meanwhile,
  work in  inverse scattering theory had found solutions to the
  classical Yang-Baxter equation associated to Lie algebras \cite{BD82}
  and  efforts to quantize these were just bearing fruit
  \cite{Drinfeld83,Jimbo85}.  This remarkable confluence resulted
  in the construction of quantum groups and associated link
  invariants, which we sketch in the next two sections.

\section{Ribbon Hopf Algebras and Tangle Invariants}

The modern view of algebraic topology is as being about certain
functors from the geometric category of topological spaces and
homotopy classes of maps, to the algebraic category of groups and
homomorphisms.  The tremendous success of this program might suggest
considering other geometric categories and hoping to find functors from
them into algebraic categories.  This  approach may be
taken as the guiding philosophy behind the invariants discussed in this
article.  In fact, the principal difference between these modern
functorial invariants and the classical ones is that here the
topological spaces themselves form the morphisms.

In this section we  construct the appropriate geometric category, the
category of framed tangles \cite{Yetter88,Turaev89}, and show that the
Kauffman bracket arises out of a functor from this category to the
category of vector spaces.   We  then show how to construct such a
functor from the representation theory of a ribbon Hopf algebra, which
we will define.  In the next section, we  construct the underlying
ribbon Hopf algebra, and discuss how to construct a large family of
similar ribbon Hopf algebras.  Despite the first two paragraphs, this whole
section
 should be accessible to someone knowing little or
no category theory.

A tangle is the image of a smooth embedding of a union of circles and
intervals into the cylinder $D \times I$, where $D$ is the unit disk
in $\CC$.
The intersection of a tangle with the boundary of the  cylinder is
required to be transverse, to  lie in $X \times(\{0\} \cup
\{1\})$, where $X$ is
the $x$-axis in $D$, and to be exactly the image of the endpoints of
the intervals.  Tangles are considered up to smooth isotopy of the
cylinder leaving $X\times \{0\}$ and $X \times\{1\}$ invariant.
Oriented and framed tangles are defined
by analogy, but we require the framing to point in the positive $x$
direction at the top and bottom.  We define $\dom(T)$ to be the number
of intersection points with
 $X\times \{0\}$, and $\codom(T)$ to be the number of intersection
 points with  $X \times\{1\}$ (for oriented tangles, we must also keep
 track of the orientations of the points).
See Figure \ref{fg:tangles} for an examples.

\fig
\begin{tabular*}{\textwidth}{c@{\extracolsep{\fill}} c@{\extracolsep{\fill}}
    c@{\extracolsep{\fill}} c}

\pic{Tngl-o/3dtangle1}{50}{-25}{0}{0}&
\pic{Tngl-o/3dtangle2}{50}{-25}{0}{0}&
\pic{Tngl-o/3dtanglecompose}{50}{-25}{0}{0}&
\pic{Tngl-o/3dtensor}{50}{-25}{0}{0}\\
$T_1$ & $T_2$  & $T_1 \comp T_2$ & $T_1 \tnsr T_2$
\end{tabular*}
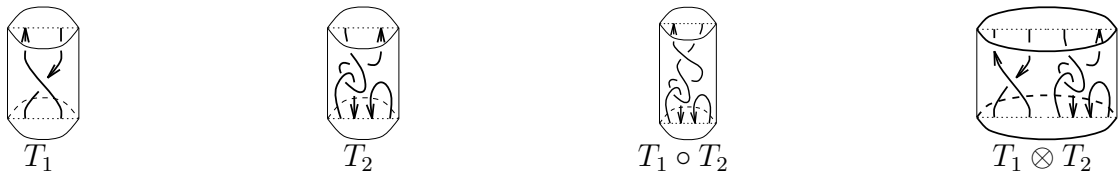
\efig{Composition and tensor product of tangles}{fg:tangles}

Tangles allow two different multiplications.  The composition of
tangles $T_1 T_2$ is the tangle formed by putting $T_1$ on top
of $T_2$, isotoping so that the boundaries match up
smoothly. This is only well-defined if $\dom(T_1)=\codom(T_2)$.  The
tensor product $T_1 \tnsr T_2$ is formed by putting them next to each
other, and treating them as a single tangle.  Both multiplications are
associative.  See Figure \ref{fg:tangles} for examples.

In the language of category theory, let $\frak{T}$ be the category whose
objects
are nonnegative integers, and whose morphisms from $n$ to $m$ are the
tangles $T$ with $\dom(T)=n$ and $\codom(T)=m$.  Composition is of
course composition of tangles.  The tensor product makes this category
into what is known as a  monoidal category.  The analogy to keep in
mind  is the category $\frak{V}$, whose objects are
finite-dimensional vector spaces and whose morphisms are linear maps.
Here the monoidal structure is given by tensor product of linear maps
and vector spaces.

In fact every tangle (hence every link, which is
a tangle from the object $0$ to itself) can be constructed out of a
handful of generating tangles.  Further, we can say exactly which such
combinations give the same tangle.  This gives a completely algebraic
presentation of this monoidal category.  The proof of this theorem is
not difficult:  It involves projections of tangles, and
paths between projections.

\thm{th:presentation}  \cite{Turaev88,Turaev89,FY89}  Every
unoriented, framed tangle is the composition of tensor
  products of the five tangles $\chr{Tngl-u/tovercross}$,
  $\chr{Tngl-u/tundercross}$,
  $\chr{Tngl-u/tcup}$, $\chr{Tngl-u/tcap}$,
  $\chr{Tngl-u/tvert}$, and two such products correspond to the same
  tangle if and only if
  they can be connected by a sequence of the following moves

\vspace{0pt}
\begin{minipage}[t]{2.5in}
\begin{Romanlist}
\vspace{0pt}
\item $
\pic{Tngl-u/lmoveI}{36}{-14}{2}{0} =
\pic{Tngl-u/rmoveI}{36}{-14}{2}{0}$
\vspace{3pt}
\item $
\pic{Tngl-u/lmoveII}{24}{-10}{2}{0}=
\pic{Tngl-u/rmoveII}{24}{-10}{2}{0}$
\vspace{3pt}
\item $
\pic{Tngl-u/lmoveIII}{36}{-14}{2}{0}=
\pic{Tngl-u/rmoveIII}{36}{-14}{2}{0}$
\end{Romanlist}\end{minipage}
\begin{minipage}[t]{3in}
\vspace{-10pt}

\begin{Romanlist} \setcounter{Numeral}{3}
\vspace{7pt}
\item $
\pic{Tngl-u/lmoveIV}{24}{-10}{2}{0}=
\pic{Tngl-u/tvert}{24}{-10}{2}{0}=
\pic{Tngl-u/rmoveIV}{24}{-10}{2}{0}$
\vspace{3pt}
\item $
\pic{Tngl-u/lmoveV}{24}{-10}{2}{0}=
\pic{Tngl-u/rmoveV}{24}{-10}{2}{0}$
\vspace{3pt}
\item $
\pic{Tngl-u/lmoveVI}{24}{-10}{2}{0}=
\pic{Tngl-u/tangle}{20}{-8}{2}{0}=
\pic{Tngl-u/rmoveVI}{24}{-10}{2}{0}$
\vspace{3pt}
\item $
\pic{Tngl-u/lmoveVII}{24}{-10}{2}{0}=
\pic{Tngl-u/rmoveVII}{24}{-10}{2}{0}$
\end{Romanlist}
\end{minipage}
\vspace{4pt}
\newline
where $T$ and $S$ are arbitrary tangles.  The same is true of oriented,
framed tangles, if each generator and each relation above is written
with every possible consistent orientation.  The same is also true
in either case for unframed tangles, if one adds in relation I that
both equal the identity tangle.
\qed \ethm

Much of this structure is quite natural algebraically.
Moves VI and VII are just restatements of axioms of a monoidal category.
More interestingly,  the morphism $\chr{Tngl-u/tovercross}$, which permutes
two objects, looks
 like the canonical vector space map from $V\tnsr W$ to $W \tnsr V$,
except that its square is not one.  The vector space map  is an
example of a {\em symmetry}, a concept  arising naturally in category
theory (see Section 5).  Instead, $\chr{Tngl-u/tovercross}$ and
$\chr{Tngl-u/tundercross}$ give the category of tangles a
{\em braiding}, a weakening of symmetry that comes up in studying
2-categories.  In fact, the subcategory of $\frak{T}$ generated by
$\chr{Tngl-u/tovercross}$, $\chr{Tngl-u/tundercross}$ and
$\chr{Tngl-u/tvert}$ is the free braided, monoidal
category generated by one object.  This subcategory is none other than
the braid group mentioned in Section 1 (actually, the braid groupoid:
The endomorphisms of the object $n$ form the mapping class group of
the plane with $n$ points removed).  The morphisms $\chr{Tngl-u/tcap}$ and
$\chr{Tngl-u/tcup}$ are  like the canonical vector space map from $V\tnsr
V^*$ to $\CC$ and its dual, and in particular give the category a
rigidity or dual structure, also much-studied in category theory.
Leaving off Move I entirely (this corresponds to Kauffman's `regular
isotopy' of links) and orienting the tangles we have the free rigid,
braided, monoidal category on
one object.  There are hints that Move I is also natural  in the
context of $n$-categories \cite{FY89,JS93,KV??}.

We are looking for a monoidal functor from $\frak{T}$ to
$\frak{V}$: That is, an assignment of linear maps to the five
generating tangles, such that the seven relations of the previous
theorem hold true, with composition interpreted as composition of
linear maps and tensor product interpreted as  tensor
product of linear maps. Call such a map a {\em tangle functor}. In particular,
a
tangle functor would send links to linear
maps from $\CC$ to itself ($\CC$, like the object $0$, is the
identity object for tensor multiplication), which are just  complex
numbers. Thus we get a numerical link invariant.

Hypothesizing that the Kauffman bracket arises from  a tangle functor, we
see that Equations (\ref{eq:Kskein}) and (\ref{eq:Kunknot}) are just two
more equations that our five linear maps must satisfy.  With patience
and the
additional hint that all of these equations can be satisfied by
operators on a two-dimensional space, you would undoubtedly come up
with something like the following \cite{Turaev89}: Let $V$ be $\CC^2$,
and define maps
\begin{align*}
\bracket{\chr{Tngl-u/tovercross}}: V \tnsr V \to V \tnsr V\quad & \quad
\bracket{\chr{Tngl-u/tundercross}}: V \tnsr V \to V \tnsr V\\
\bracket{\chr{Tngl-u/tcap}}:  V \tnsr V \to \CC\quad
\bracket{\chr{Tngl-u/tvert}}: V  &\to V \quad
\bracket{\chr{Tngl-u/tcup}}: \CC \to V \tnsr V
\end{align*}
by the following, writing the basis for $V\tnsr V$ as $(1,0)\tnsr
(1,0)$, $(1,0)\tnsr (0,1)$, $(0,1) \tnsr (1,0)$, and $(0,1)\tnsr (0,1)$:
\begin{align}\label{eq:Ktf1}
\bracket{\chr{Tngl-u/tovercross}} =
    \begin{bmatrix}
    A & 0 & 0 & 0\\
    0& 0 & A^{-1} & 0\\
    0 & A^{-1} & A - A^{-3} &  0\\
    0 & 0 & 0 & A
    \end{bmatrix} \,\, & \quad
\bracket{\chr{Tngl-u/tundercross}}=
    \begin{bmatrix}
    A^{-1} & 0 & 0 & 0\\
    0 & A^{-1} - A^3 & A & 0\\
    0& A & 0 & 0\\
    0 & 0 & 0 & A^{-1}
    \end{bmatrix}\\
\pagebreak[3] \label{eq:Ktf2}
\bracket{\chr{Tngl-u/tcap}}=
    \begin{bmatrix}
      0 & A & -A^{-1} & 0
     \end{bmatrix}
\quad
\bracket{\chr{Tngl-u/tvert}}&=
      \begin{bmatrix}
        1 & 0\\
        0 & 1
        \end{bmatrix}
\quad
\bracket{\chr{Tngl-u/tcup}}=
      \begin{bmatrix}
        0 &
        -A &
        A^{-1} &
        0
        \end{bmatrix}^T.
\end{align}

The reader can check all the relations with nothing but elementary
linear algebra (Moves III and V are best checked by applying Equation
(\ref{eq:Kskein}) first).  In particular this confirms that the
Kauffman bracket is a framed link invariant. Using the above the Jones
polynomial can also be written as a tangle functor.

We now give a set of structures and axioms  on an algebra ${\cal A}$
which guarantee that the category whose objects are finite-dimensional
representations and whose morphisms are intertwiners (maps between
representations commuting with the action of ${\cal A}$) is naturally
the range of such a functor.  Let us proceed heuristically.

 Since sequences of points are
sent to a tensor product of representations, the algebra must act on
such tensor products.  Likewise, it must act on $\CC$, the image of
the object $0$.  This requires a {\em bialgebra\/}.  That is, an
algebra ${\cal A}$, with homomorphisms $\Delta:{\cal A}  \to {\cal A}
\tnsr {\cal A}$ (the coproduct) and $\epsilon:{\cal A} \to \CC$ (the
counit) satisfying
\begin{align*}
(\Delta \tnsr 1)\Delta&= (1\tnsr \Delta)\Delta\\
(\epsilon \tnsr 1)\Delta&= 1= (1\tnsr \epsilon)\Delta.
\end{align*}
The reason for the terminology is that the adjoint of $\Delta$
together with $\epsilon$ give a multiplication and identity on the
dual space, making it an algebra.  A bialgebra ${\cal A}$
acts on $\CC$ by $\epsilon: {\cal A} \to \CC= \End(\CC)$, and if
$\rho:{\cal A} \to \End(V)$ and $\tau:{\cal A} \to \End(W)$ are
representations, then $(\rho
\tnsr \tau)\Delta:{\cal A} \to \End(V) \tnsr \End(W)$ is a
representation on $V\tnsr W$.

Now restricting  attention to oriented tangles, recall that
$\chr{Tngl-u/tcap}$ and $\chr{Tngl-u/tcup}$, with various orientations,
were analogous to duality in vector
spaces.  That is,  if a boundary point with one orientation
corresponds to a representation $V$, a point with the other
orientation should correspond to $V^*$, and for example
$\chr{Tngl-o/forcap}$ should
correspond to the canonical map from $V^*
\tnsr V$ to $\CC$. In order for $V^*$ to be a
representation,  our bialgebra ${\cal A}$ must be a {\em Hopf algebra}:
It should have an antihomomorphism $S:{\cal A} \to {\cal A}$
(antipode) satisfying
$$m(S\tnsr 1)\Delta=m(1\tnsr S)\Delta=1_{{\cal A}}\epsilon,$$
where $m$ is the multiplication map.  If $\rho:{\cal A} \to \End(V)$ is
a representation,
then the dual representation $\rho^*:{\cal A} \to \End(V^*)$ is defined by
$\rho^*(a)
v^*=v^*\comp \rho(S(a))$.  The antipode axiom assures that there is a
canonical intertwiner from $V\tnsr V^*$ to the trivial representation.
Hopf algebras have a long and distinguished history
\cite{Sweedler69,LR87} outside the scope of this article.

The last piece of information is the braiding $\chr{Tngl-u/tovercross}$, which
should be an intertwiner from $V\tnsr W$ to $W \tnsr V$.  It should
not be the flip map,  $\sigma_{VW}:v\tnsr w \mapsto w\tnsr v$, because then
we would have $\chr{Tngl-u/tovercross}=\chr{Tngl-u/tundercross}$ and our tangle
functor would be very uninteresting!  The flip map will be an
intertwiner whenever ${\cal A}$ is cocommutative, that is, when
all $\Delta(a)$ are symmetric.  Thus we will expect interesting link
invariants only from noncocommutative Hopf algebras.  More precisely, if
$\Delta'(a)= \sigma_{{\cal A}{\cal A}}(\Delta(a))$, then we expect $\Delta\neq
\Delta'$.
They shouldn't be too different, though, since we want
$\chr{Tngl-u/tovercross}$   to
behave somewhat like the flip map.  It turns out the right notion
is for $\Delta$ and $\Delta'$ to be connected by an inner
automorphism.  That is, define
a {\em quasitriangular Hopf algebra} ${\cal A}$ \cite{Drinfeld85} to be a Hopf
algebra
${\cal A}$ together
with an invertible element $R \in {\cal A}\tnsr {\cal A}$, satisfying
\begin{align*}
\Delta'(a)&=R\Delta(a)R^{-1}\\
(\Delta\tnsr 1)R=R_{13}R_{23}\quad & \quad
(1\tnsr \Delta)R=R_{13}R_{12}
\end{align*}
where if $R=\sum_i a_i \tnsr b_i$, then $R_{13}= \sum_i a_i \tnsr 1
\tnsr b_i$, $R_{23}=\sum_i 1 \tnsr a_i \tnsr b_i$,  etc.  From these
 equations one easily gets the Yang-Baxter equation
\eq{eq:YB}
R_{12}R_{13}R_{23}=R_{23}R_{13}R_{12}
\end{equation}
which will ultimately give us Move III.

Notice cocommutative Hopf
algebras are always trivially quasitrangular, with $R=1\tnsr 1$.
Drinfeld has shown \cite[4.2D]{CP94}\cite{Drinfeld87} that one can
combine any Hopf algebra
${\cal A}$ with its dual to get a larger Hopf algebra, the `quantum
double,' which is quasitriangular (the ones we are interested in
are almost of this form).  Thus there are lots of them around.  They
also have a rich and interesting algebraic structure.  For example, if
we define $u= \sum_i S(b_i)a_i$, where $R=\sum_i a_i \tensor b_i$,
then
$$S^2(a)=uau^{-1}.$$
This implies for example that $V$ is isomorphic as a representation to
its double dual
$V^{**}$, though not by the obvious map.

The representation theory of a quasitriangular Hopf algebra will be a
rigid, braided, monoidal category.  To get
 Move I requires a {\em ribbon Hopf algebra\/} \cite{RT90}. The ribbon
 structure is not as compelling algebraically as the Hopf and
 quasitriangular structure, perhaps because we do not understand Move
 I categorically as well as the others. A ribbon Hopf
algebra is a quasitriangular Hopf algebra with a grouplike element $G$
(i.e., $ \Delta(G)= G\tnsr G$ and $\epsilon(G)=1$, hence $S(G)=G^{-1}$)
satisfying
\begin{align*}
G^{-1}uG^{-1}&=S(u)\\
GaG^{-1}&=S^2(a).
\end{align*}
 Given a representation
$(\rho_V,V)$ and an element $x\in \End(V)$, define the quantum trace
$\qtr$ by $\qtr_V(x)=
\tr(\rho_V(G)x)$.  The second condition above implies that, identifying
$\End(V)$ with $V \tnsr V^*$, $\qtr$ is the unique to scaling
intertwiner from this representation to the trivial one.  That $G$ is
grouplike means that the quantum trace is additive on direct sums and
multiplicative on tensor products.  Also important is the quantum
dimension of $V$, $\qdim(V)= \qtr_V(1)$.

A simple example is  called for.  Let $G$ be a finite group, and
$\CC G$ be the group algebra.  Then the maps
$$\Delta(g)= g\tnsr g\quad
S(g)=g^{-1}\quad
\epsilon(g)=1 \quad (\text{for}\,g \in G)$$
extend by linearity to make $\CC G$ a cocommutative Hopf algebra.  This
explains
the term grouplike. Being cocommutative, it is a ribbon Hopf algebra with
$$R=1\tnsr 1 \qquad
G=1.$$
While these are certainly interesting Hopf algebras, their quasitriangular and
ribbon
structures are trivial, and we should not expect information about links
from them.  For that we must wait until the next section.

In putting this together, we see that we have gotten more than we were
originally looking for.  Specifically, given a ribbon Hopf algebra
${\cal A}$, define a {\em labeled tangle} to be an oriented tangle with a
finite-dimensional representation of ${\cal A}$ assigned to each
component.  Thus the domain and codomain are now sequences of oriented
points labeled by representations.  Assign to such a sequence a
representation of ${\cal A}$ by
$${\cal F}((V_1,\varepsilon_1),(V_2,\varepsilon_2), \ldots,
(V_n,\varepsilon_n) )= \bigotimes_{i=1}^n V_i^{\varepsilon_i},$$
where $V_i$ is a representation, $\varepsilon_i$
is $+$ or $-$ indicating the orientation, and $V^+_i=V_i$,
$V^-_i=V_i^*$.
\thm{th:tanglefunctor}\cite{RT90}
For every ribbon Hopf algebra ${\cal A}$ there is a monoidal functor from
the category of labeled oriented framed tangles to that of
representations and intertwiners of ${\cal A}$.  That is, there is
a map ${\cal F}$ which assigns to each such tangle $T$ an intertwiner
${\cal F}(T) : {\cal F}(\dom(T)) \to {\cal F}(\codom(T))$, satisfying
${\cal F}(T_1T_2)= {\cal F}(T_1){\cal F}(T_2)$ and ${\cal F}(T_1\tnsr
T_2)= {\cal F}(T_1)\tnsr {\cal F}(T_2)$, ${\cal F}$ of the identity
  tangle is the identity operator, and ${\cal F}(\emptyset)=1$.  It is
  determined by  its values on
  generating tangles, which are shown in the following chart.  Here
  $G_i$ and $R_{ij}$ are the actions of $G$ and $R$ on $V_i$ and $V_i
  \tnsr V_j$ respectively, $\sigma_{ij}$ is the flip map, and
  $v_\alpha$ is a basis of $V_i$, $v_\alpha^*$ its dual basis of
  $V_i^*$.  Mirror images of the crossings shown get sent to their
  inverses.
\ethm

\begin{pf}
We have only to check invariance under Moves I-VII.  Moves II, IV, VI
and VII are easy to check, Move III is the Yang-Baxter Equation
(\ref{eq:YB}).  Moves I and V require more involved calculations.
\end{pf}

\begin{tabular}{|| l | l | l || l | l | l ||} \hline
${\cal F}$ of & takes & to &${\cal F}$ of & takes & to \\ \hline
$\rule{0pt}{14pt}\chr{Ltngl-o/backcap}$& $x^*\tnsr x$  & $x^*(x)$
&$\chr{Ltngl-o/backcup}$ & $c$ & $c\sum_\alpha v_\alpha \tnsr
v_\alpha^*$ \\ \hline
$\rule{0pt}{14pt}\chr{Ltngl-o/forcap}$& $x \tnsr x^* $ &
$x^*(G_i(x)) $ &$\chr{Ltngl-o/forcup}$ & $c $& $c\sum_\alpha
v^*_\alpha \tnsr
G_i^{-1}(v_\alpha) $ \\ \hline
$\rule{0pt}{14pt}\chr{Ltngl-o/ppcross}$& $x\tnsr y $ &
$\sigma_{ji}R_{ij} (x \tnsr y)$ & $\chr{Ltngl-o/mmcross}$&
$x^*\tnsr y^*$ & $\sigma_{j^*i^*}R_{i^*j^*}(x^* \tnsr y^*)$
\\ \hline
$\rule{0pt}{14pt}\chr{Ltngl-o/pmcross}$& $x\tnsr y^* $ &
$\sigma_{j^*i}R_{ij^*} (x \tnsr y) $ & $\chr{Ltngl-o/mpcross}$&
$x^*\tnsr y  $ &
$\sigma_{ji^*}R_{i^*j} (x^* \tnsr y) $
\\ \hline
\end{tabular}

\begin{remark}\mbox{}
\begin{itemize}
\item There is a kind of converse of this result, a Tannaka-Krein type
theorem:  Every functor from $\frak{T}$ to $\frak{V}$ arises in this
way from some ribbon Hopf algebra \cite{Majid89,Sawinnotes}.
\item We have focused on oriented tangle functors, although currently
  our only interesting example is unoriented.  It turns out that a
  self-dual representation of a ribbon Hopf algebra gives (almost) an
  unoriented tangle functor.  As is shown in the next section,  the
  Kauffman bracket functor essentially arises from a Hopf
  algebra representation corresponding to the fundamental
  representation of $sl_2$, which is self-dual.
\item The use of tensor products in discussing $\Delta$ and $R$ is
  subtle in the case of most interest,when ${\cal A}$ is infinite-dimensional.
  All equations make sense if we
   require only that they hold when represented on an
  arbitrary finite-dimensional representation.  This involves no loss
  of information, because in the end we are only interested in these
  representations \cite{Sawinnotes}.
\end{itemize}
\end{remark}

We need to beef up our invariant in a fairly trivial fashion for
future sections.  Define a {\em ribbon graph} exactly as a tangle,
except we allow {\em coupons}, squares with a definite top, bottom and
orientation, with strands allowed to intersect the coupons at distinct
points on the top and bottom, as in Figure \ref{fg:rgraphs}.  Theorem
\ref{th:presentation} is still true of this larger category if we add
the coupon as a generator and add relations VIII (and its mirror
image) and IX in Figure
\ref{fg:rgraphs}, understood to apply for coupons with any number of
incoming and outgoing strands.  Given a ribbon Hopf algebra ${\cal A}$,
a {\em labeled ribbon graph} is one where the strands are labeled by
representations of ${\cal A}$ as before, and any coupon with incoming
strands labeled $V_1,\ldots,V_n$ and outgoing labeled by
$W_1,\ldots,W_m$ is labeled by an intertwiner from $V_1 \tnsr \cdots
\tnsr V_n $ to $W_1 \tnsr \cdots \tnsr W_m$, with representations
replaced by their duals if the corresponding strand is oriented down
when it intersects the coupon.  It is an easy matter to show that
${\cal F}$ extends to a functor from labeled ribbon graphs to the
representation theory of ${\cal A}$, with ${\cal F}$ of a coupon labeled
by $f$ being simply $f$.  The following summarizes important facts
about ${\cal F}$.

\fig
\begin{tabular}{c|c@{\hspace{10pt}}c|c}
An Example \hspace{4pt}& \multicolumn{2}{c}{Rules} & Labeling\\
\rule{0pt}{5pt}&&&\\
\pic{Coupon/example}{40}{-20}{0}{0} &
\hspace{4pt}$\pic{Coupon/lmoveVIII}{40}{-20}{0}{0}=\pic{Coupon/rmoveVIII}{40}{-20}{0}{0}$
&$\pic{Coupon/lmoveIX}{40}{-20}{0}{0}=
\pic{Coupon/rmoveIX}{40}{-20}{0}{0}$\hspace{4pt}
&\pic{Coupon/label}{40}{-20}{0}{0}\\
\rule{0pt}{2pt}&&&\\
 & VIII  & IX & \hspace{4pt}$f:V \to U \tnsr W^*$
\end{tabular}
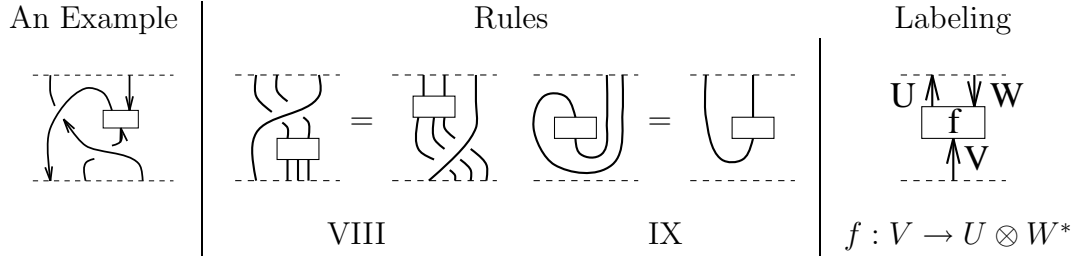
\efig{Coupons in ribbon graphs}{fg:rgraphs}

\prop{pr:facts}\cite{RT90}
\begin{alphalist}
\item ${\cal F}$ of a tangle with a component labeled by a direct sum of
  representations is the direct sum of ${\cal F}$ of the tangle labeled by each
  summand.
\item ${\cal F}$ of a tangle with a component labeled by a given
  representation is ${\cal F}$ of the same tangle with the component
  given the reverse orientation and labeled by the dual
  representation.
\item ${\cal F}$ of a tangle with one component labeled by the trivial
    representation is ${\cal F}$ of that tangle with that component
    deleted.
\item ${\cal F}$ of a tangle with a component labeled by a tensor
  product of representations is ${\cal F}$ of that tangle with the
  component replaced by two parallel components, following the
  framing, labeled by the tensor factors.
\item ${\cal F}$ of a tangle which can be separated by a sphere into a
  link and a subtangle is ${\cal F}$ of the link times ${\cal F}$ of the
  subtangle.
\item Let $T$
  be a tangle with a component labeled by an irreducible
  representation $V$, and let $L$ be a link with a component
  labeled by $V$.  If $T'$ is the tangle formed by cutting both
  of these components and gluing them together along the cuts
  \rom{(}consistent with the orientations\rom{)} then ${\cal F}(T')=
  {\cal F}(T){\cal F}(L)/\qdim(V)$.
\item ${\cal F}(\chr{Coupon/qtrf})=\qtr(f)$.\qed
\end{alphalist}
\renewcommand{\qed}{\mbox{}}
\eprop

\section{Quantum Groups}

Up to this point, this article has been close to being self-contained:
All that has been left out are computations, and with a few more
precise statements and a sprinkling of hints, the industrious reader
could probably reproduce them.  I hope also that it has been
 well-motivated: That the idea of a functor from the tangle
category to  vector spaces is a natural thing to look for, and
that it might reasonably lead one to consider ribbon Hopf algebras.
The same will not be true of this section.

It is certainly natural to argue, as we do, that Lie algebras give a
nice set of  ribbon Hopf algebras which offer trivial link
information because they are cocommutative.
This suggests deforming them within the set of all ribbon
Hopf algebras to get  interesting  link
invariants.  If we could argue that there was a unique deformation in
the space of ribbon Hopf algebras, or that some sort of quantization
described it geometrically, such a construction would  seem very natural.
Unfortunately, these approaches give only a rough framework in which
to proceed, and ultimately one is forced to write down a guess for
generators and relations and prove by hand that they give a ribbon
Hopf algebra.  For this reason we will eschew the deep and significant
mathematics in these approaches
\cite[\S1-3,6.1,6.2]{CP94}\cite{FT87,Drinfeld83,Drinfeld87}, and
present the algebras as if by oracle.  Likewise, there is a fair
amount of interesting technical machinery to prove the basic facts of
these algebras
\cite[\S6.4,8.1-8.3,10.1]{CP94}\cite{Lusztig93,Lusztig88,Rosso90b,Jimbo85,KR88,KR90}:
We will merely sketch the easy case
$sl_2$ and state the broad results
in general.  Ultimately, from our point of view these quantum groups are a
means
for extracting geometric information (link invariants) out of
quintessentially geometric objects (Lie groups), and it is certainly
to be hoped that in the near future we will have a geometric
understanding of how they arise.

Recall from the previous section that the group algebra of a finite
group is a ribbon Hopf algebra, although the most information it can
give about a link is the number of components.  It is a little subtle
deciding how to define the analogue of the group algebra for a Lie
group, but it is clear that morally the same should be true.  The best
surrogate for this object turns out to be the universal enveloping
algebra of the Lie algebra, $U(\frak{g})$.  For the group $SL(2)$,
this is the algebra $U(sl_2)$ generated by $\{x,y,h\}$, with relations
\begin{align*}
[h,x]&=2x\\
[h,y]&=-2y\\
[x,y]&=h
\end{align*}
where $[a,b]=ab-ba$.  Its Hopf algebra structure is determined by
$$
\Delta(a)=a \tnsr 1 + 1 \tnsr a\qquad
\epsilon(a)=0\qquad
S(a)=-a
$$
where $a$ is in $\{x,y,h\}$.  Its ribbon and quasitriangular
structure is trivial, $R=1 \tnsr 1$ and $G=1$.  The universal
enveloping algebras of the other simple Lie
algebras have similar, if more complicated presentations
\cite{Humphreys72}.   Given such a Lie algebra $\frak{g}$, there is
associated an algebra $U_s(\frak{g})$, the quantum universal enveloping
algebra of $\frak{g}$, depending on a nonzero complex parameter $s$.  It will
prove to be a ribbon Hopf algebra, and in an appropriate sense will
approach $U(\frak{g})$ as a ribbon Hopf algebra as $s \to 1$. These
quantized universal enveloping algebras are collectively called quantum groups,
though the term is sometimes used more generally to refer to ribbon
Hopf algebras, quasitriangular Hopf algebras, or even all Hopf
algebras.

$U_s(sl_2)$ is the algebra generated by $\{x,y,h\}$
subject to the relations
\begin{align}
[h,x]&=2x \nonumber\\
[h,y]&=-2y\\
[x,y]&=(s^{2h}-s^{-2h})/(s^{2}-s^{-2}).\nonumber
\end{align}
There are a number of artifices for interpreting the last equation,
none entirely satisfactory (the casual reader is encouraged to ignore
the problem).  Ours is as follows:  Let $U(sl_2)'$ be the
  completion of $U(sl_2)$ in the topology of convergence on every
  finite-dimensional representation, and let the closure
  of the subalgebra generated by $h$ be  ${\cal H}$. Notice  the
  right side of the last equation
is in ${\cal H}$ for all $s \neq 0$, because every finite-dimensional
representation is
spanned by eigenvectors of $h$ with integer eigenvalues.  Consider the
algebra spanned by free products of $x$, $y$ and elements of ${\cal H}$,
with the topology inherited from ${\cal H}$.  The above equations
generate an ideal in this algebra, and the quotient by the closure of
that ideal is the algebra $U_s(sl_2)$.

This algebra is really just $U(sl_2)$ in disguise.  Specifically, if
we use $\bar{x}$, $\bar{y}$, $\bar{h}$ temporarily to denote the
generators of $U(sl_2)$ and recall that the Casimir element is
$C=(\bar{h}+1)^2/4 + \bar{y}\bar{x},$ then the map $\varphi:U_s(sl_2) \to
U(sl_2)'$ given by
\begin{align*}
\varphi(h)=\bar{h}\quad & \quad
\varphi(y)=\bar{y}\\
\varphi(x)=4(s^{4\sqrt{C}}+s^{-4\sqrt{C}}-
s^{2\bar{h}-2}+s^{2-2\bar{h}})&/
((\bar{h}-1)^2 - 4C)(s^2-s^{-2})^2 \bar{x}
\end{align*}
is a homomorphism for all $s\neq 0$.  Further, this homomorphism is
1-1 and has
dense range if $s$ is not a root of unity (essentially because the
coefficient of $\bar{x}$ is invertible) \cite[\S4.6]{CP94}\cite{Jimbo85}.
Thus up to completions
 $U_s(sl_2)$ and $U(sl_2)$ are isomorphic as algebras.
So $U_s(sl_2)$ is semisimple, its finite-dimensional
representations are in one to one correspondence with that of
$U(sl_2)$, and the actions of $h$ and $\bar{h}$ on these representations are
the same.
In particular, they have the same formal characters.

The Hopf algebra structure is
given by
\begin{equation}
\begin{array}{rlcrl}
\Delta(x)&= x\tnsr s^h + s^{-h} \tnsr x & \qquad & S(x)&=-s^2x\\
\Delta(y)&= y\tnsr s^h + s^{-h} \tnsr y & \qquad & S(y)&=-s^{-2}y
\end{array}
\end{equation}
and the rest as in the $U(sl_2)$ case. The key observation is that
these maps commute with the operator $[h,\cdot]$.  From this and the fact that
the
formal characters are the same as for $U(sl_2)$ it is easy to see that
duals of representations and decomposition of tensor products into
irreducible representations are the same (away from a root of unity)
as in the classical case.

 The quasitriangular structure
is quite a bit trickier.  Define quantum integers
$$
[n]_q= (s^{2n}-s^{-2n})/(s^2 - s^{-2}).
$$
Then
\begin{equation}R= s^{(h\tnsr h)} \sum_{n=0}^\infty
  \frac{s^{n(n+1)}(1-s^{-4})^n}{[n]_q[n-1]_q
\cdots[1]_q } x^n \tnsr y^n\qquad
G=s^{2h}.\end{equation}
 The infinite sum is well defined in the usual topology, because only
 finitely many terms are nonzero on any given pair of representations.

If you believe  that all the axioms of a ribbon Hopf algebra are
satisfied by the algebra above, you should be wondering what the
tangle functor looks like.  First we must understand the
representation theory explicitly.  For this recall that
$sl_2$ has an $n$-dimensional representation for each $n>0$, with
basis $\{v_i\}_{i=1}^n$ such that
\begin{align*}
hv_i&=(n-2i+1)v_i\\
xv_i&=(i-1)v_{i-1}\\
yv_i&=(n-i)v_{i+1}
\end{align*}
and that the action of $h$ is independent of $s$ (by the isomorphism).
It is easy to reconstruct the representation $V_n$ of $U_s(sl_2)$ with
basis $\{v_i\}_{i=1}^n$ such that
\begin{align}
hv_i&=(n-2i+1)v_i \nonumber\\
xv_i&=[i-1]_qv_{i-1} \label{eq:repdef}\\
yv_i&=[n-i]_qv_{i+1}.\nonumber
\end{align}

It is now easy to compute the action of $R$ on $V_2 \tnsr V_2$, and
composing with the flip map gives
$$\chr{Tngl-o/ppcross} \to \begin{bmatrix} s & 0 & 0 & 0\\
0 & 0 & s^{-1} & 0\\
0 & s^{-1} & s-s^{-3} & 0\\
0 & 0 & 0 & s
\end{bmatrix}.$$
Identifying $V$ with $V^*$ by the intertwiner $\alpha:v_1 \mapsto
sv_2^*$, $\alpha:v_2 \mapsto -s^{-1}v_1^*$ gives
$$
\chr{Tngl-o/forcap} = \begin{bmatrix} 0 & -s & s^{-1} & 0
\end{bmatrix}\qquad
\chr{Tngl-o/backcap} =\begin{bmatrix}0 & s & -s^{-1} & 0
\end{bmatrix}
$$
with the other generating tangles determined by these.  Notice this is
{\em almost\/} an unoriented tangle functor---in fact the Kauffman
bracket functor, Equations (\ref{eq:Ktf1}) and
(\ref{eq:Ktf2}), with  $A=s$---except for an annoying minus sign in
$\chr{Tngl-o/forcap}$. The link invariant will actually be the
Kauffman bracket times $(-1)^w$, where $w$ is the total winding number
of the link projection ($(-1)^w$ is a framed, unoriented link
invariant!).  In general any irreducible self-dual representation
gives, on choosing such
an $\alpha$, either an unoriented functor or a functor which differs
from one by a minus sign in $\chr{Tngl-o/forcap}$.  Kirillov and
Reshetikhin \cite{KR88} nicely fix this problem for $sl_2$ by
constructing $U_q(sl_2)$ as a deformation of $U(sl_2)$ with a
nonstandard quasitriangular (really  triangular) structure.

We have
only described the functor with  everything labeled by the
two-dimensional representation, but we have really constructed much
more than that.  There is enough information above to
write down the full theory for all representations, but  no neat
skein-theory description has been found for the corresponding invariant.

It is now clearer where the skein relation in Equation (\ref{eq:Kskein})
is coming from.  $V_2 \tnsr V_2$ breaks up as the sum of two
irreducible representations, the trivial one and $V_3$, and thus the space
of intertwiners on it is two-dimensional.  The two fragments
pictured on the right side of Equation (\ref{eq:Kskein}), interpreted
as tangles,
correspond to the identity intertwiner and a multiple of the projection
onto the trivial subrepresentation.  They span the space of intertwiners,
which includes $\chr{Tngl-u/tundercross}$.

As for other Lie algebras, the story is essentially the same.  A
similar, but more complicated presentation of $U_s(\frak{g})$ can be
given, and shown in the same sense to be isomorphic to $U(\frak{g})$
(though in general the isomorphism does not admit explicit formulae as
in the $sl_2$ case).  The ribbon Hopf structure can all be written down
explicitly, and the representation theory is `the same' as for the
original Lie algebra except at roots of unity
\cite[\S8.3,10.1]{CP94}.  The link invariants
coming from the fundamental representation of $sl_n$ are particular
values of the two-variable HOMFLY polynomial.  The polynomials for
$B_n$, $C_n$, and $D_n$ at the fundamental representation are all
special values of the two-variable
Kauffman invariant \cite{Turaev89,Reshetikhin88}.  Both of
these polynomials can be computed by
skein-theoretic algorithms similar to that for the Jones polynomial.
There is also a skein-theory algorithm for the $G_2$ invariant at the
fundamental representation involving trivalent graphs
\cite{Kuperberg94}.

\begin{remark} A word on notation.  It is customary to speak of the
  quantum groups as $U_q(\frak{g})$, where $q$ is, depending on the
  author, $s^2$ \cite{CP94,Kassel94,Lusztig93,Jimbo85} or $s^4$
  \cite{Reshetikhin88,Rosso90,Drinfeld87,KM91}.  Of course, this requires
choosing
  a second
  or fourth root of $q$ to write the $R$ matrix.  It is also common to
  call this variable $t$, or to call the variable in the Jones
  polynomial $q$, although these are inverses of each
  other with the usual conventions.  Other minor variations exist on
  the definition of $U_q(\frak{g})$.  We have tried to remain internally
  consistent and consistent with well-established conventions, such as
  the skein relations for the Jones polynomial.
\end{remark}

\section{Modular Hopf Algebras}

What happens to the quantum group at a root of unity?  The answer is,
the algebra ceases to be isomorphic to the unquantized algebra, and in
fact ceases to be semisimple.  The ribbon
Hopf structure and some of the representation theory survive this
collapse however, and emerge with subtle properties that allow the link
invariant to be extended to a $3$-manifold invariant and a topological
quantum field theory.  This is the
subject of the next four sections.  We begin with a precise look at
the situation for $sl_2$.

When $s$ is a root of unity, the representations of $U_s(sl_2)$ given
in the previous section are still representations, but  may not all be
irreducible.  It is not hard to see that irreducibility is equivalent
to $v_1$ being the unique $h$-eigenvector on which $x$ acts as $0$.
Thus if we let $l$ be the least natural number such that $s^{4l}=1$
(so that $s$ is a primitive $l$th or $2l$th root of unity for $l$
odd, or a primitive $4l$th root of unity), we have $[l]_q=0$, $[n]_q
\neq 0$ for $n<l$, and thus $V_n$ is still irreducible for $n \leq l$,
but is not for larger $n$, by Equation (\ref{eq:repdef}).  Recalling
that the quantum dimension of a representation is
 the  trace of
$G$ in that representation,
notice that the quantum dimension
of $V_n$ is $[n]_q$, and thus is nonzero for $n<l$ and equal to zero
for $n=l$.  It turns out that irreducible representations with quantum
dimension zero have an important property: If one takes any
representation in the ideal generated by them, i.e., a direct sum of
tensor products of these with other representations, then any
intertwiner from this representation to itself will have quantum trace
zero.  But this implies that any link labeled by such a representation
has invariant zero.  It is easy to check that the tensor product $V_n
\tnsr V_m$ for $n,m<l$ is a sum of representations $V_k$ for $k<l$ plus
a representation in this ideal.  Since for purposes of knot theory the
representations in this ideal are irrelevant, it makes sense to throw
this trivial-trace representation out and define the {\em truncated
  tensor product} to be the
rest of the tensor product.  With this new tensor product, direct sums
of representations $V_n$ for $n<l$ still form a rigid, braided, tensor
category, and give a tangle functor with truncated tensor product.
But now there are only finitely
many essentially different labels.

The same situation applies to a general quantum group, at least in
many cases.  If $s$ is a primitive $2l$th root of unity for $l$ odd
(and also prime to $3$ in the $G_2$ case), then the
 representation of highest weight $\lambda$ is still  an irreducible
representation with the correct character when
$\bracket{\lambda+\rho,\theta} \leq l$, where $\theta$ is the highest
root (the unique long root in the Weyl chamber), $\rho$ is half the
sum of the positive roots, and
$\bracket{\cdot,\cdot}$ is as in \cite{Humphreys72}.
This representation has nonzero quantum dimension exactly when the
inequality is strict.  Further, if the decomposition of the unquantized
tensor product of two representations in the smaller set
lies entirely in the larger, the tensor product in this case
decomposes the same way.  From this it follows that the tensor product
of any two representations in the smaller set is a
sum of such representations plus a trivial-trace representation.
Again, this gives a `truncated' tangle
functor \cite{Andersen92}.  Of course, this is only interesting if
more than the trivial representation is included, so we usually
restrict to $l$ such that $l>\bracket{\rho,\theta}$.  The quantity
$l-\bracket{\rho,\theta}$ is the $k$ occurring in the Chern-Simons
action in the introduction.

The other useful  property of quantum groups at roots of unity has to
do with the value of the invariant on the Hopf link, B in Figure
\ref{fg:linkexamples}.  Let $H_{i,j}$ be the value of ${\cal F}$ on the
Hopf link labeled by representations $\lambda_i$ and $\lambda_j$.  For
$U_s(sl_2)$, a straightforward calculation shows that $H_{i,j}=
[ij]_q$, where of course $\lambda_i=i$.  At a root of unity, if we
restrict to our preferred representations, these numbers form a matrix,
and another calculation shows that this matrix is nonsingular when $s$ is a
primitive $4l$th root of unity or a primitive $2l$th root with $l$
odd \cite{RT91,BHMV92}.  Turaev and Wenzl have shown \cite{TW93}, by
relating the
question to work of Kac and Petersen \cite{KP84}, that the same is
true of this matrix for an arbitrary quantum group when $l$ satisfies
the same restrictions as in the previous paragraph.

This information is encapsulated in the definition of a modular Hopf
algebra (see \cite{RT91}, in this simplified form \cite{Turaev92,TW93}).

\begin{definition}
A {\em modular Hopf algebra}  is a ribbon Hopf algebra ${\cal A}$, together
with a finite collection of irreducible representations $
\lambda_1, \ldots, \lambda_k$, including the trivial representation,
 closed under duals, and having  $\qdim(\lambda_i)\neq
0$, satisfying
\begin{itemize}
\item $\lambda_i \tnsr \lambda_j = \bigoplus_{m=1}^k (N_{i,j}^m
  \lambda_m) \oplus \eta_{i,j}$, where $N_{i,j}^m$ is a multiplicity
  and $\eta_{i,j}$ is a representation on which $\qtr$ is zero on all
  intertwiners.
\item The matrix $H_{i,j}$, for $1 \leq i,j \leq k$ is
  nonsingular.
\end{itemize}
\end{definition}

\begin{remark}
That the quantum dimension of each $\lambda_i$ is nonzero actually
follows from  the second condition of the
definition.
\end{remark}

We will now derive the key facts about a modular Hopf algebra that we
will need later on.  Consider the truncated representation algebra
$\bold{R}$ of our
modular Hopf algebra. That is, the commutative algebra over $\CC$
spanned by
$\lambda_1, \ldots, \lambda_k$ with $\lambda_i \lambda_j = \sum_m
N_{i,j}^m \lambda_m$.

We can extend $\qdim$ by linearity to a homomorphism from
$\bold{R}$ to $\CC$.  More generally, if $T$ is a ribbon tangle which is
labeled
except for one closed component $x$ without coupons, and $a= \sum_i c_i
\lambda_i \in \bold{R}$ we can define $T$ with $x$ labeled by $a$
to be the formal sum $\sum_i c_i T_{\lambda_i}$, where $T_{\lambda_i}$
is $T$ with $x$ labeled by $\lambda_i$.  Then define the value of the
invariant
$${\cal F}(T_a) \defequals \sum_{i=1}^k c_i {\cal F}(T_{\lambda_i}).$$
This is consistent with the additivity of the invariant on direct
sums.

It is  natural to ask if $\bold{R}$ is semisimple:  That is, is it
spanned by minimal idempotents?   In fact the dual basis to
$\{\lambda_i\}$ with respect to the nondegenerate, symmetric pairing
given by $\{H_{i,j}\}$ consists of such minimal idempotents.  Minimal
idempotents correspond to homomorphisms to $\CC$, and the idempotent
dual to the trival representation corresponds to the homomorphism
$\qdim(\cdot)$.  This minimal idempotent, which is the source of the
$3$-manifold invariant, is computed explicitly in the following
proposition.

\prop{pr:Omega} \cite{Sawinnotes} Let
$$\omega= \sum_{i=1}^k \qdim(\lambda_i)\lambda_i.$$
\begin{alphalist}
\item $a\omega = \qdim(a)\omega$ for all $a \in \bold{R}$.
\item the Hopf link labeled by $\lambda_i$ and $\omega$ has nonzero
  ${\cal F}$ if and only if $\lambda_i$ is trivial \rom{(}in particular
  $\qdim(\omega) \neq 0$\rom{)}.
\qed
\end{alphalist}
\eprop
\section{Axiomatic Topological Quantum Field Theory}

The fundamental idea of the first three sections was to get a detailed
combinatorial description of an interesting geometric category
($\frak{T}$) and use it to find functors to categories of vector
spaces and linear maps.  We  do the same thing  with the
cobordism category, and call it a topological quantum field theory, or
TQFT.  Now we  obtain numerical invariants not of (framed) links,
but of (biframed) $3$-manifolds.  The work is remarkably parallel,
the key difference being that this category is not presented as
straightforwardly, but
indirectly, in terms of links.  Thus instead of TQFT's being
associated to a specific kind of algebra, we  construct them from
a specific kind of link invariant: one arising from a modular Hopf
algebra.

Let  $\Sigma_1$ and $\Sigma_2$ be smooth oriented  $d-1$-manifolds.
A dimension $d$ cobordism $\frak{m}$  with domain $\Sigma_1$ and
codomain $\Sigma_2$
is up to diffeomorphism a triple
$(M,f_1,f_2)$, where $M$ is an oriented smooth $d$-dimensional manifold with
boundary, and  $f_1$ and $f_2$ are orientation-preserving endomorphisms
from $\Sigma_1^*$ and $\Sigma_2$ respectively to $\partial M$.
Here $\Sigma_1^*$ is $\Sigma_1$ with the opposite orientation, and we
require of $f_1$ and $f_2$ that
 $\partial M$ be the disjoint union of their ranges.  By ``up to
diffeomorphism'' we mean that two triples $(M,f_1,f_2)$ and $(M',f'_1,f'_2)$
are the same morphism if there is a diffeomorphism $F:M_1 \rightarrow
M_2$ with $f'_1 = F f_1$ and $f'_2 = F f_2.$  Thus it is a manifold
with parameterized boundary, but with some of the boundary considered
incoming and some considered outgoing, so that it can be viewed
 as acting like a function.  More precisely, think of the $\Sigma_i$
as being analogous to vector spaces and cobordisms as being analogous
to linear maps.

In this vein, the analogue of composition of linear maps is gluing of
cobordisms.  That is, if $\frak{m} = (M,f_1,f_2)$
and $\frak{m}' = (M',f'_1,f'_2)$, define  $\frak{m}'\frak{m} = (M'
\cup_{f_1'f_2^{-1}} M, f_1,f_2')$, where $M' \cup_{f_1' f_2^{-1}} M$
is the manifold formed by identifying points in $\partial M'$ with
$\partial M$ via the orientation reversing map $f_1' f_2^{-1}$.
This composition is associative, and its identity is clearly
$$1_{\Sigma} = (\Sigma \times I, \mbox{id}_0, \mbox{id}_1),$$
where $I$ is the unit interval and
$\mbox{id}_0:\Sigma^* \rightarrow \Sigma \times \{ 0 \}$ and
$\mbox{id}_1:\Sigma \rightarrow \Sigma \times \{ 1 \}$ are the identity
maps.

Tensor product of vector spaces is analogous to disjoint union
$\Sigma_1 \cup \Sigma_2$; $\CC$, the identity for tensor product,
corresponds to $0$, the empty $d-1$-manifold.
Corresponding to tensor product of linear maps we have the associative
operation $\frak{m} \cup \frak{m}' = (M \cup M', f_1 \cup f'_1,f_2
\cup f'_2)$.  The empty $d$-manifold $\emptyset$ acts like the
number $1$ as the identity for this product.

Of course $V\tnsr W$ and $W \tnsr V$ are the same vector space, in the
sense that there is a canonical map between them $\sigma_{VW}: v
\tnsr w \mapsto w \tnsr v$.  Likewise there is a canonical cobordism from
$\Sigma_1 \cup \Sigma_2$ to $\Sigma_2 \cup \Sigma_1$, namely
the cobordism $\frak{c}_{\Sigma_1,\Sigma_2} = ((\Sigma_1 \cup
\Sigma_2) \times I, f_0,f_1),$ with
$$f_0:\Sigma_1^* \cup \Sigma_2^* \rightarrow (\Sigma_1 \cup \Sigma_2)
  \times \{ 0 \} \quad f_1:\Sigma_2 \cup \Sigma_1 \rightarrow (\Sigma_1 \cup
\Sigma_2) \times
 \{ 1 \}$$
being the identity and the order-reversing map respectively.

No mathematician, and certainly no category theorist, would let such a
good analogy go unnamed.  The definition of cobordisms is
designed to give the morphisms of a category, which we'll call
$\frak{C}$, and all the additional structures we discussed are just to
say that $\frak{C}$ and $\frak{V}$ are
symmetric, monoidal or {\em tensor\/} categories
(there are some axioms to check, which are trivial in both cases. See
\cite{MacLane71}).

The point of this is
\begin{definition} \cite{Atiyah89,Atiyah90b}
  A {\em $d$-dimensional axiomatic topological quantum field theory,}
  or {\em TQFT}, is a functor ${\cal Z}$ of tensor categories
  from $\frak{C}$ to $\frak{V}$.  That is, a map ${\cal Z}$ which sends
  each oriented $d-1$-manifold to a finite-dimensional vector space
  ${\cal Z}(\Sigma)$, and each cobordism $\frak{m}$ from $\Sigma_1$ to
  $\Sigma_2$ to a linear map ${\cal Z}(\frak{m}):{\cal Z}(\Sigma_1) \to
  {\cal Z}(\Sigma_2)$, such that ${\cal Z}(1_\Sigma)=1_{{\cal Z}(\Sigma)}$,
  ${\cal Z}(\Sigma_1 \cup \Sigma_2)= {\cal Z}(\Sigma_1) \tnsr
  {\cal Z}(\Sigma_2)$, ${\cal Z}(0)=\CC$, ${\cal Z}(\frak{m}_1 \cup
  \frak{m}_2)= {\cal Z}(\frak{m}_1) \tnsr {\cal Z}(\frak{m}_2)$,
  ${\cal Z}(\emptyset)=1$, and ${\cal Z}(\frak{c})$ is the flip map.
\end{definition}

Notice that a closed $d$-manifold is a cobordism from the empty
$d-1$-manifold to itself, and thus gets sent to an operator on
 $\Bbb{C}$, or in other words a number.  So one thing a
TQFT gives is a numerical $d$-manifold invariant.  There is a
straightforward algorithm for computing it from a Morse function: If
you can compute the operators assigned to cobordisms for each type of
singularity, you can multiply them together to get the invariant of
the closed manifold.  Also, it follows  that
the mapping class group of any $d-1$-manifold acts on the vector space
associated to that manifold.  Thus these contain a lot of topological
information \cite[\S III]{Turaev94}.

The definition of a TQFT is modeled on what we expect a
topological quantum field
theory to yield.  For the passage from the
physics to these axioms,  see \cite{Atiyah90b} or \cite{Axelrod91}.

Let us now fix $d=3$.  We will not actually
construct a functor from $\frak{C}$ to $\frak{V}$, but from the
category $\frak{F}\frak{C}$ of {\em biframed} (also called $2$-framed)
cobordisms to $\frak{V}$.
The reason for this lies buried in the subtleties of the Chern-Simons
path integral, from which these examples spring, where a choice of
biframing is necessary to regularize certain integrals.  The casual
reader may safely ignore the technicalities of biframings, and view
them as playing a role very similar to that of framings on links for
link invariants.

A {\em biframing} on a closed $3$-manifold is a trivialization of $TM \oplus
TM$, up to isotopy.  Atiyah argues \cite{Atiyah90} that a given
trivialization extends
to a $4$-manifold bounded by that $3$-manifold exactly when the
$4$-manifold has a certain signature, which is a complete invariant of
the biframing.  A biframing on a manifold with boundary
$\Sigma$ restricts to a  map
from $T\Sigma \oplus T\Sigma $ to $\RR^6$. Call a
{\em biframed $2$-manifold\/} one equipped with such a map.  A biframed
cobordism is a $3$-cobordism together with  a trivialization up to
isotopy fixing the boundary, and its domain and codomain are thus
biframed $2$-manifolds.  Biframed cobordisms form a compact closed
category $\frak{F}\frak{C}$ exactly as before, with the observation
that if $\Sigma$ is a
biframed $2$-manifold, $\Sigma^*$ inherits its biframing, and $1_\Sigma$,
 and $ \frak{c}_{\Sigma,\Gamma}$
 have canonical biframings.

We can make this category simpler to work with in a number of ways.
First, two biframed $2$-manifolds of the same genus are easily seen to
be connected by an invertible biframed cobordism.  Thus their vector spaces are
isomorphic, and we only need to pick a vector space for each one (this
reduction is what category theorists call `skeletonization'.  The
category of tangles is really a skeletonization  of a larger category
with a geometrically more natural definition).
The fact that $\frak{F}\frak{C}$ has a duality structure similar to
  that of tangles (view $\Sigma \times I$ as a morphism from $\Sigma^*
  \cup \Sigma$ to $0$) means we really do not have to distinguish
  between domain and codomain.
Finally, the value of the
functor on disconnected cobordisms is clearly determined by its value
on connected ones.

To state the simplified version, for each $g$ choose a representative
biframed genus $g$ surface
$\Sigma_g$, and a biframed cobordism $\frak{d}_g: \emptyset \to
\Sigma_g \cup \Sigma_g$ whose underlying manifold is $\Sigma_g \times
I $ and which is symmetric in the sense that
$\frak{c}_{\Sigma_g,\Sigma_g}\frak{d}_g = \frak{d}_g$.  Also let
$\frak{e}_g$ be a biframed cobordism such that $(\frak{e}_g \tnsr 1)(1 \tnsr
\frak{d}_g)
 = 1 = (1 \tnsr \frak{e}_g)(\frak{d}_g \tnsr 1)$.

\thm{th:TQFTclass} \cite{Sawinnotes} For each $g$, let
${\cal Z}(\Sigma_g)$ be a finite-dimensional vector space;  for each
biframed cobordism $\frak{m}:\bigcup_{i=1}^n \Sigma_{g_i} \to \emptyset$,
let ${\cal Z}(\frak{m}): \bigotimes_{i=1}^n
{\cal Z}(\Sigma_{g_i}) \to \Bbb{C}$ be a linear functional;  and for each
nonempty, closed, biframed $3$-manifold $M$, let ${\cal Z}(M)$ be in
$\Bbb{C}$.  Then ${\cal Z}$ extends uniquely to
a biframed TQFT, if and only if the following hold:
\begin{alphalist}
\item {\em Nondegeneracy:\/}  ${\cal Z}(\frak{e}_g)$ is a symmetric
  nondegenerate pairing on
  ${\cal Z}(\Sigma_g)$.

\item {\em Symmetry:\/} If $\frak{m}$ has domain $\bigcup_{i=1}^n
  \Sigma_{g_i}$, and $\frak{m}'$ with domain $\bigcup_{i=1}^n
  \Sigma_{g_{\sigma(i)}}$ for some permutation $\sigma$ is the same
  manifold as $\frak{m}$ with the same parameterization written in a
  different order, then ${\cal Z}(\frak{m})= {\cal Z}(\frak{m}')
  P_\sigma$, where $P_\sigma$ is the map on the appropriate tensor
  product of vector spaces which permutes the tensor factors.

\item {\em Sewing:\/} Suppose $\frak{m}$ has domain $\bigcup_{i=1}^m
  \Sigma_{g_i}$ and $\frak{n}$ has domain $\bigcup_{i=1}^n
  \Sigma_{g'_i}$ with $g_m=g'_1$.  We can
  form the biframed manifold
  $\frak{m}\cup_s \frak{n}$, with domain $\bigcup_{i=1}^{m-1}
  \Sigma_{g_i} \cup \bigcup_{i=2}^n \Sigma_{g'_i}$, by composing with
  $\frak{d}_g$ along these boundary components (Figure
  \ref{fg:sewmendex}).  Then
$${\cal Z}(\frak{m}\cup_m \frak{n})=
{\cal Z}(\frak{m})\tnsr {\cal Z}(\frak{n})\circ \alpha$$
where $\alpha$ is the canonical map sending $v_1 \tnsr \cdots \tnsr
v_{m-1} \tnsr w_2 \tnsr
    \cdots \tnsr w_n$ to $\sum_i v_1 \tnsr \cdots v_{m-1} \tnsr a_i
    \tnsr b_i \tnsr w_2 \tnsr
    \cdots \tnsr w_n$, with  $a_i$ and
    $b_i$ dual bases of ${\cal Z}(\Sigma_{g_m})$ with respect to the
    pairing ${\cal Z}(\frak{e}_g)$.

\item {\em Mending:\/} Suppose $\frak{m}$ has domain $\bigcup_{i=1}^n
  \Sigma_{g_i}$, with $g_1=g_2$.  We can form $\frak{m}_m$ by
  composing with
  $\frak{d}_g$ along the first two boundary components (Figure
  \ref{fg:sewmendex}).  Then
$${\cal Z}(\frak{m}_m)=
{\cal Z}(\frak{m})\circ \alpha$$
where $\alpha$ is
the canonical map as above.\qed
\end{alphalist}  \ethm

\fig
\begin{tabular*}{5.4in}{c@{\extracolsep{\fill}}c}
$\pic{Misc/lsewex}{50}{-30}{0}{0} \rightarrow \pic{Misc/rsewex}{30}{-10}{0}{0}$
& $\pic{Misc/lmendex}{50}{-30}{0}{0} \rightarrow
\pic{Misc/rmendex}{30}{-10}{0}{0} $\\
Sewing & Mending
\end{tabular*}
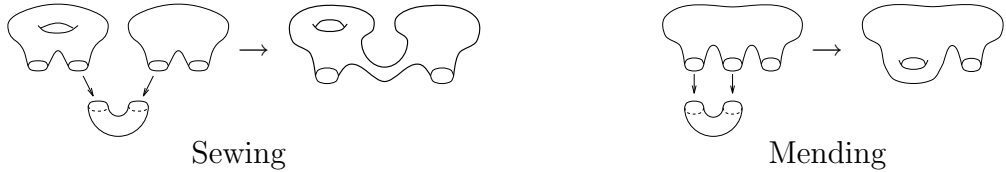
\efig{Examples of sewing and mending}{fg:sewmendex}

\begin{remark} \mbox{}
\begin{itemize}
\item There is one obvious piece of structure on $\frak{C}$ that we are
leaving off:  Reversing the orientation on a cobordism from $\Sigma$
to $\Sigma'$ gives a cobordism from $\Sigma'$ to $\Sigma$.  This is
analogous to Hilbert space adjoint.  It thus seems natural to replace
$\frak{V}$ with the category of finite-dimensional Hilbert spaces and
ask that this structure too be preserved.  This is called a {\em
  unitary} TQFT, and is the more natural one from the point of view of
physics.  It turns out that, while we will construct a TQFT for any
quantum group and any primitive $2l$th root of unity (with the
limitations on $l$ given in Section 4) it is only unitary for
$s=e^{\pi i/l}$.   This is
exactly the value of $s$ corresponding to Witten's geometric construction.
\item Of course, an ordinary framing on a closed $3$-manifold
  determines a biframing, so in the end we will have an invariant of
  framed $3$-manifolds as well.
\end{itemize}
\end{remark}
\section{Surgery and Cobordism}

We would like to have an algebraic description of the category
$\frak{F}\frak{C}$ of framed cobordisms in terms of generators and
relations.  The idea is to  define a TQFT by
the image of the generating morphisms, and prove
that it is a TQFT by confirming that it preserves the
relations.  A  natural choice is Morse theory.  Here the
generators are the morphisms attaching handles and the relations are
those given by Cerf theory \cite{Cerf70}.  This is the approach taken
in \cite{Walker??}, but it
involves a lot of detail checking, because there are a lot of ways to
attach a handle.

Surgery on links meets our present need much better.  Since surgery
describes 3-manifolds in terms of links, and we are constructing
TQFT's out of link invariants, we will find many of the details fall
into place.  Of course, surgery does not apparently handle biframing, or
manifolds with boundary, so this section is devoted to extending it to
these situations.  We begin with a review of (integer) surgery.

Consider a framed unoriented link $L$ in $S^3$. Let $T$ be a tubular
neighborhood of $L$.  Each component of $\partial T$ has a natural
meridian, which bounds a disk in $T$.  It also has a natural
longitude,  the pushoff of $L$ in the direction of the framing.
Remove each component of $L$ and glue it back in by a map sending the
meridian to the longitude and the longitude to minus the meridian (you
must choose orientations for these curves, but the result does not
depend on this choice).  The manifold one obtains is called $M_L$, the
result of surgery on $L$.  For example, surgery on the $0$-framed
unknot is easily seen to be $S^1\times S^2$, and  a little more
work shows that surgery on a $\pm 1$ framed unknot
gives $S^3$ again.

It turns out, by a theorem of Lickorish \cite{Lickorish62}, that every
closed, compact, connected, oriented $3$-manifold admits such a
presentation.  Further,
by a theorem of Kirby \cite{Kirby78}, two links present the same
$3$-manifold if and only if they can be related by a sequence of the
moves in Figure \ref{fg:Kirby} and their mirror images, with any
number of strands understood to pass through the pictured unknot.
Actually, Kirby's original version involved two ``nonlocal'' moves:
Fenn and Rourke \cite{FR79} reduced them to this move.

\fig
$$\pic{Kirby/lmovekirby}{55}{-25}{0}{0}=\pic{Kirby/rmovekirby}{55}{-25}{0}{0}$$
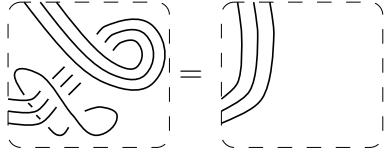
\efig{The Kirby move}{fg:Kirby}

Kirby's proof actually takes place in four dimensions.  He views $S^3$
as bounding $B^4$, and interprets $L$ as defining a $4$-manifold by
attaching $2$-handles to tubular neighborhoods of the components.  The
resulting $4$-manifold will have boundary $M_L$.  Thus the more
sophisticated view of surgery is that it results in a $3$-manifold
together with a choice of $4$-manifold for it to bound.

This is perfect for us, because this is just what we saw in the
previous section determines a biframing.  Of course  we
need  a refinement of the Kirby move which preserves the signature
of the $4$-manifold but is still  powerful enough to connect all links
giving the same signature.  This task at first sounds daunting, but is
in fact quite simple.  The key observation is that, since each
component of $L$ corresponds to adding a handle to the $4$-manifold,
the second relative cohomology of the $4$-manifold $M$ has
 a basis element for each handle, supported in that handle.
Further, the intersection form of two basis elements is
 the linking number of the corresponding components (the
self-intersection is the self-linking number, determined by the
framing).  Thus the signature of
the $4$-manifold is exactly the signature of the matrix of linking
numbers of $L$!

\fig
$$\mbox{I.} \quad
\pic{Kirby/framed-unknots}{45}{-20}{0}{0}=\pic{Dfrag-u/dashbox}{45}{-20}{0}{0}
\qquad \qquad
\mbox{II.} \quad
\pic{Kirby/lmovekirby}{55}{-25}{0}{0}=\pic{Kirby/rmovefkirby}{55}{-25}{0}{0}$$
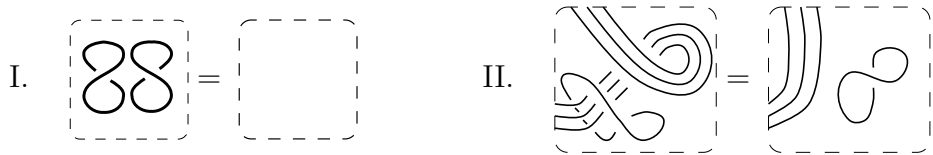
\efig{The biframed Kirby moves}{fg:FKirby}

An easy calculation shows that the Kirby move changes the signature of
the linking matrix by $1$.  Thus the two {\em biframed Kirby moves\/}
pictured in Figure~\ref{fg:FKirby}, each a composition of two Kirby
moves, do not change the signature.  On the other hand, if two links
represent the same $3$-manifold and have the same signature, one can
convert the sequence of ordinary Kirby moves connecting them to a
sequence of biframed Kirby moves.
\thm{th:FKirby}\cite{Sawinnotes} Closed, oriented, biframed
$3$-manifolds are in one to
one correspondence with equivalence classes of biframed unoriented links
in $S^3$ modulo moves I and II in Figure~\ref{fg:FKirby} and their
mirror images. \qed
\ethm

The cobordism category is not much harder.
For each surface $\Sigma_g$ choose a biframed handlebody which it bounds,
$H_g$.  To represent pictorially an embedding of $H_g$ into $S^3$,
 choose  a set of generators for the fundamental group of the
handlebody, $\{c_i\}$ for $1\leq i \leq g$.  Represent an embedding of
$H_g$ into $S^3$ as in
Figure~\ref{fg:embedding}.The embedding is recovered by thickening
the framed graph to a handlebody and identifying $H_g$ with the thickening by a
map
which sends  $c_i$ to the boundary above the $i$th circle on the
graph.  Choose $\frak{d}_g$ so that the identification of $\Sigma_g$
with $\Sigma^*_g$ sends each $c_i$ to a meridian intersecting it and
vice versa, as illustrated for genus two in Figure \ref{fg:embedding}.

\fig\begin{tabular}{c|c}

$$\pic{Hndlbd/handlebody-embedding}{80}{-30}{0}{0}=
\pic{Hndlbd/embedding-presentation}{80}{-30}{0}{0}$$
&\pic{Coupon/mapping-cylinder}{60}{-25}{0}{0}\\
\rule{0pt}{16pt}Presenting handlebody embeddings
\hspace{40pt}&\hspace{40pt} Boundaries
identified by $\frak{d}_g$
\end{tabular}
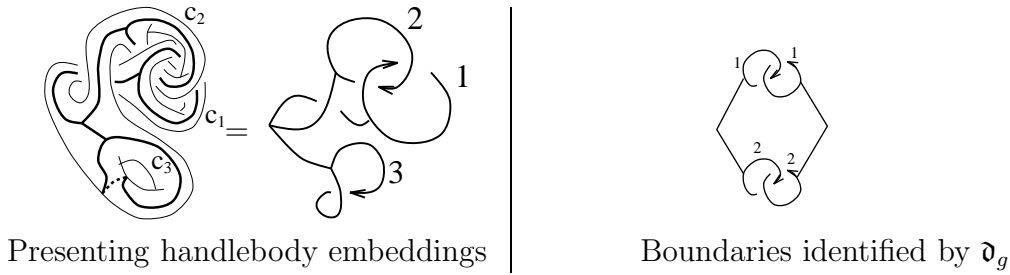
\efig{A pictorial description of the cobordism category}{fg:embedding}

Now if $\frak{m}: \bigcup_{i=1}^n \Sigma_{g_i} \to \bold{\emptyset}$
is a connected cobordism, consider the closed manifold $M$ formed by gluing
$\bigcup_{i=1}^n H_{g_i}$ to $\frak{m}$ along the boundary by the
parametrization, and consider the resulting embedding $N$ of
$\bigcup_{i=1}^n H_{g_i}$ into $M$.  This pair determines $\frak{m}$
uniquely.  $M$ can be presented as $M_L$ for some link $L$, and the
embedding $N$ can be isotoped so as  not to  intersect the embedded
tori coming from surgery on $L$. so $N$ corresponds to an embedding of
$\bigcup_{i=1}^n H_{g_i}$ into $S^3$ which does not intersect $L$.
Thus $\frak{m}$ is determined by a pair $(N,L)$, where $L$ is an
unoriented biframed link and $N$ is an embedding of handlebodies into
$S^3$ which does not intersect $L$.  Specifically, $\frak{m}$ is
obtained by doing surgery on $L$, removing the image of the interiors
of the handlebodies, and parametrizing the resulting boundary by $N$
restricted to the boundary. We call $(N,L)$ a presentation of
$\frak{m}$, and in general refer to the  cobordism presented by
$(N,L)$ as $[N,L]$.

In \cite{Roberts??b} a set of Kirby-type moves is given for manifolds
with boundary.  From this it follows that two presentations give the
same cobordism if and only if they can be connected by isotopy and a sequence
of
 Kirby moves (Figure~\ref{fg:Kirby}) and their mirror image.  The framed
unknot in Figure~\ref{fg:Kirby}  represents a component of $L$ and the strands
passing
through it represent pieces of $L$ or $N$.  We thus have a purely
combinatorial description of cobordisms.

In fact, $(N,L)$ determines a biframed cobordism, and
Figure~\ref{fg:FKirby} connects all equivalent presentations of a
biframed cobordism.
 Specifically, having chosen a
biframing on  $H_g$, biframings on
$[N,L]$ are in one-to-one correspondence with biframings on $M_L$, which
are preserved by the biframed Kirby moves.

Our last goal  is to give a presentation of a cobordism
obtained by sewing or mending, in terms of the presentation of the
manifold(s) to be sewn or mended.  The following
proposition is not hard to show from the biframed Kirby moves.

\prop{pr:standard-position} \cite{Sawinnotes} Any biframed cobordism can be
written as $[N,L]$,
where $N$ is in {\em standard position\/} in the sense that it can be
projected with no self-crossings, as illustrated in
Figure \ref{fg:sandm}.
\eprop

Thus we may  assume the presentations to be sewn or mended are in
standard form.  It is more clear and convenient to describe the algorithm
pictorially, as in Figures~\ref{fg:sandm}.  The
corresponding picture applies for any genus and any number of strands
passing through the handles.
\fig \begin{tabular}{c|c}
$\begin{array}{cc}\pic{Coupon/lsewing}{55}{-25}{0}{0}&\Rightarrow
\pic{Coupon/rsewing}{55}{-25}{0}{0}\\
\rule{0pt}{24pt} A \hspace{60pt} B &  A \cup_s B
\end{array}$ \hspace{20pt}&\hspace{20pt}
$\begin{array}{cc}
\pic{Coupon/lmending}{70}{-30}{0}{0}&\Rightarrow
\pic{Coupon/rmending}{70}{-30}{0}{0}\\
\rule{0pt}{16pt}A & A_m
\end{array}$
\end{tabular}
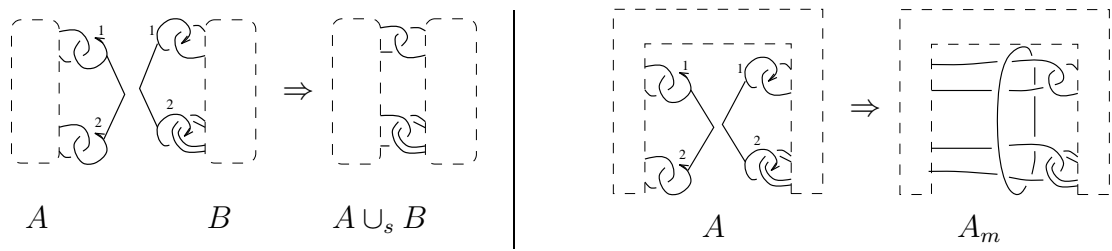
\efig{Sewing and mending}{fg:sandm}

The algorithm for sewing follows  straightforwardly from the
definition, but the algorithm for mending needs a little
justification.  We first note that, to form the mending $[N,L]_g$, we
might just as well form the mending $[N,\emptyset]_m$ and then do
surgery on the image of $L$ in this new manifold.  But now because $N$
is in standard position, $[N,\emptyset]_m$ is simply $S^1 \times S^2$,
and the picture we get is exactly the right side of Figure~\ref{fg:sandm}.

\section{Constructing The TQFT}

We are now ready to construct a TQFT.  For this we
need to associate vector spaces to surfaces.  We  first associate
excessively large vector spaces to each $\Sigma_g$, spanned by labeled
ribbon graphs in $H_g$.  Then we  give
a pairing on this corresponding to $\frak{e}_g$, and quotient by
the null space to force the pairing to be nondegenerate.  The value of
the invariant on cobordisms is then fairly obvious.

Given a biframed handlebody $H_g$, a closed labeled ribbon graph in $H_g$ is
defined  as in Section 2, except that the graph is
embedded in $H_g$, and equivalence is by a
positive diffeomorphisms fixing the boundary.

Let $V_g$ be the vector space of all formal linear combinations of
closed labeled ribbon graphs in $H_g$.  For each presentation $(N,L)$,
where $N$ is an embedding of $\bigcup_{j=1}^n H_{g_i}$ into $S^3$, we
define a map $f_{(N,L)}:\bigotimes_{j=1}^n V_{g_j} \to \CC$, as
follows.  If $h_j$ is a closed labeled ribbon graph in $H_{g_j}$ for
$1\leq j \leq n$, embed $\bigcup_{j=1}^n h_j$ into $S^3$ via $N$, and
label each component of $L$ by $\Omega=
\qdim^{-1/2}(\omega)\omega$, with $\omega$ as in
Proposition~\ref{pr:Omega}.
Let $K=\qdim(\Omega)$, and define $f_{(N,L)}(\bigotimes h_j)= K^{-1}$ times
${\cal F}$ of the resulting graph.
\prop{pr:invariance} \mbox{}
\begin{alphalist}
\item If Moves I and II of Figure~\ref{fg:FKirby} represent closed labeled
ribbon graphs in
  $S^3$ with the framed unknots pictured labeled by $\Omega$, the
  values of ${\cal F}$ on  both sides of the equal signs agree.
\item The function $f_{(N,L)}$ depends only on $[N,L]$, and not on the
  presentation.  Thus we will speak henceforth of $f_{[N,L]}$.
\end{alphalist} \eprop

\begin{pf}
\begin{alphalist}
\item
  For Move II, use Proposition \ref{pr:facts} and Figure
  \ref{fg:sequence}.  The first step is by parts (a) and (d), the second
  by part (f), the third by Proposition (\ref{pr:Omega}a), and the last
  by (a) and (d) again.

For Move I, apply Move II to the left  side to get a Hopf link with
both components labeled by $\Omega$, one having a $+1$ framing.  By
Proposition (\ref{pr:Omega}b), ${\cal F}$ of this is $\qdim
(\Omega)\qdim^{-1/2}(\omega)=1$.
\item This follows immediately from (a). \qed
\end{alphalist}
\renewcommand{\qed}{}
\end{pf}
\fig
\begin{align*}
{\cal F}\left(\pic{Kirby/kirby-invariance1}{40}{-18}{-2}{2}\right) &= \sum_i
N_{\lambda_i} {\cal
F}\left(\pic{Kirby/kirby-invariance2}{40}{-18}{-2}{2}\right)=
\sum_i N_{\lambda_i}
{\cal F}\left(\pic{Kirby/kirby-invariance3}{40}{-18}{-2}{2}
\right)/\qdim(\lambda_i)\\
&= \sum_i N_{\lambda_i}
{\cal F}\left(\pic{Kirby/kirby-invariance4}{40}{-18}{-2}{2}\right) =
{\cal F}\left(\pic{Kirby/kirby-invariance5}{40}{-18}{-2}{2}\right)
\end{align*}
\efig{Proof of invariance under Move II}{fg:sequence}
 In
particular, we have a symmetric bilinear pairing $f_{\frak{e}_g} :
V_g \tnsr V_g \to \CC$.  If $N_g = \{v \in V_g: f_{\frak{e}_g}(v,w)=0
\,\,\forall w \in V_g\}$, define

\eq{eq:define-Zspace} {\cal Z}(\Sigma_g)= V_g/N_g,
\end{equation} and then $f_{\frak{e}_g}$ descends to a nondegenerate pairing on
${\cal Z}(\Sigma_g)$.
\prop{pr:well-defined} $f_{[N,L]}:\bigotimes_j V_{g_j} \to\CC$
descends to a map ${\cal Z}([N,L]):\bigotimes_j {\cal Z}(\Sigma_{g_j})
\to \CC$.
\eprop

\begin{pf*}{Sketch of Proof} This follows easily from the fact that we
  can treat $f_{[N,L]}(h_1 \tnsr \cdots \tnsr h_n)$ as
  $f_{\frak{e}_{g_j}}(h_j \tnsr k_j)$, where $k_j$ is the image of
  the other $h_i$'s and $L$ under an identification of the complement
  of the image of $H_{g_j}$ with  $H_{g_j}$.
\end{pf*}

We are now closing in on our prey.  The only question remaining is the
elusive one of sewing and mending.

\lem{lm:s-and-m}
 For any sequence of signed, labeled points, there exists a set of
  intertwiners $x_i$ and $y_i$, $1 \leq i \leq k$ such that the
  equalities hold in Figure \ref{fg:cut} for the value of
  ${\cal F}$ on any closed, labeled, ribbon graph containing  the pictures.
  \elem
\fig
$$\pic{Coupon/lcut}{50}{-25}{0}{0}=\sum_i\pic{Coupon/rcut}{50}{-25}{0}{0}
\qquad \qquad
\pic{Coupon/lcircum}{50}{-25}{0}{0}=\sum_i\pic{Coupon/mcircum}{50}{-25}{0}{0}
= K\sum_i\pic{Coupon/rcircum}{50}{-25}{0}{0}$$
\efig{Cutting ribbon graphs}{fg:cut}
\begin{pf}
The sequence of signed labeled points corresponds to some tensor
product of representations, say $W$.  Write the identity on $W$ as a
sum $\sum_j p_j$ of projections onto irreducible subrepresentations
(plus a projection onto a trivial-trace subrepresentation which does
not affect the outcome and which we ignore).  For the first,
we may interpret the fragment $T$ of the graph as a tangle, with
${\cal F}(T)$ intertwining $W$ and the trivial representation.  Thus
${\cal F}(T)p_i$ will be zero unless $p_i$ projects onto a
trivial subrepresentation.  So we get ${\cal F}(T)= {\cal F}(T) \sum_i
p_i$, with the sum over projections
onto trivial subrepresentations.  Writing each such $p_i$ as $y_i
x_i$, with $x_i$ an intertwiner from $W$ to the trivial representation
and $y_i$ an intertwiner the other way, we have the result.

For the second, it follows  from Proposition \ref{pr:Omega} that
${\cal F}$ of the tangle shown is a multiple of the projection onto the
sum of all trivial subrepresentations.  This equates
the first and second pictures.  The last is just Proposition
\ref{pr:facts}e.
\end{pf}

The following theorem first appears in \cite{RT91}, though not in the
language of TQFT's.  It was proven for $sl_2$, using only the Kauffman
bracket formalism without quantum groups, in \cite{BHMV??}.  Complete
proofs of the theorem as stated here first appear in  \cite{Walker??} and
\cite[Thm. IV.1.9]{Turaev94}.

\thm{th:TQFT} \cite{RT91}
Given a modular Hopf algebra ${\cal A}$, the ${\cal Z}$ defined in this
section extends to a TQFT.
\ethm
\begin{pf}
  We need to check axioms (a)-(d) from Theorem~\ref{th:TQFTclass}.
  Nondegeneracy we have already checked, and Symmetry is immediate
  from the construction.

For Sewing, let $[M,L]=\frak{m}$ and $[N,K]=\frak{n}$ be two
presentations in standard
form of the cobordisms mentioned in the Sewing axiom.  We may as well assume
that each has only  one boundary component, since we can  check the equality
with the maps applied to arbitrary vectors  by gluing  ribbon
graphs into each other boundary component.  Let $J$ be the
presentation of their sewing described in Section 6, and let $j$ be
the closed labeled ribbon graph obtained by labeling every component
of $J$ by $\Omega$.  Thus
$${\cal Z}(\frak{m} \cup_s \frak{n})= K^{-1}{\cal F}(j).$$
 More precisely, since $[M,L]$ is in standard form,
we can glue the complement of the image of $M$ to the underlying
manifold of $\frak{d}_g$ and identify the result with $H_g$.
This identification sends $L$ to a ribbon
graph $h_{\frak{m}}$ in $H_g$, such that ${\cal Z}(\frak{m})(h)=
\bracket{h_{\frak{m}},h}$.  Likewise we can define $h_{\frak{n}}$ so
that ${\cal Z}(\frak{n})(h)=
\bracket{h_{\frak{n}},h}$.  Now it is clear from the construction that
$j$ is $h_{\frak{m}}$ glued to $h_{\frak{n}}$ via $\frak{d}_g$, so that
$${\cal F}(j)=K\bracket{h_{\frak{m}}, h_{\frak{n}}}.$$
On the other hand
\begin{align*} \bracket{h_{\frak{m}},h_{\frak{n}}}&=\sum_{j=1}^k
  \bracket{h_{\frak{m}},a_j} \bracket{b_j,h_{\frak{n}}}\\
  &=\sum_{j=1}^k {\cal Z}(\frak{m})(a_j){\cal Z}(\frak{n})(b_j)
\end{align*}
which gives the result.

For Mending, let $[M,L]=\frak{m}$ be a  presentation in standard form
of the cobordism mentioned in this axiom, and as above assume
$\frak{m}$ has only two boundary components, both of genus $g$.  Let
$a_j$ and $b_j$ be a basis and dual basis for ${\cal Z}(\Sigma_g)$, and
choose linear combinations of closed ribbon graphs in $H_g$ to
represent them, which we'll also call $a_j$ and $b_j$.  The proof is
then contained in Figure \ref{fg:mendingproof}.  The second equality
is by Lemma \ref{lm:s-and-m}a, the third by Proposition
\ref{pr:facts}e and the seventh by Lemma \ref{lm:s-and-m}b.  The
combinations $a_j$ and $b_j$ are represented pictorially as one graph
for conservation of subscripts: It does not effect the proof.
\end{pf}
\fig
$$\sum_j {\cal Z}(\frak{m})(a_j \tnsr b_j) = K^{-1} \sum_j
{\cal F}\left(\pic{Coupon/mend1}{30}{-12}{-1}{3}\right) = K^{-1} \sum_{i,j}
{\cal F}\left(\pic{Coupon/mend2}{50}{-23}{-1}{3}\right)=$$
$$K^{-1} \sum_{i,j}
{\cal F}\left(\pic{Coupon/mend3a}{30}{-12}{-1}{-1}\right) {\cal F}\left(
\pic{Coupon/mend3b}{30}{-12}{-1}{3}\right) = K \sum_{i,j}
\left\langle\pic{Coupon/mend4a}{25}{-10}{-1}{3},a_j\right\rangle \left\langle
\pic{Coupon/mend4b}{20}{-8}{-1}{3},b_j\right\rangle=$$
$$K \sum_i \left\langle \pic{Coupon/mend4a}{25}{-10}{-1}{3},
\pic{Coupon/mend4b}{20}{-8}{-1}{3} \right\rangle = \sum_i
{\cal F}\left( \pic{Coupon/mend5}{30}{-12}{-1}{3}\right)= K^{-1} {\cal F}
\left( \pic{Coupon/mend6}{45}{-25}{-1}{3} \right)= {\cal Z}(\frak{m}_m)$$
\efig{Proof of mending invariance}{fg:mendingproof}

\begin{remark}\mbox{}
\begin{itemize}
\item These issues of framing are easy to fix for closed manifolds.  If
$[\emptyset,L]$ is the manifold presented by $L$, and $\sigma(L)$ is
the signature of the linking matrix of $L$, then a simple computation
shows that $C^\sigma(L) {\cal Z}([\emptyset,L])$ is invariant under
framing if $C= 1/u_+(\Omega)$.  This is the Reshetikhin-Turaev
invariant, except they multiply by $K$ to get ${\cal F}$
instead of ${\cal Z}$, thus making the invariant multiplicative in
connect-sums of manifolds.  Their normalization is a bit different,
because they in effect label by $\omega$
instead of $\Omega$, and thus must correct for the number of
components as well as the signature of $L$.
\item The vector spaces ${\cal Z}(\Sigma_g)$ have been defined somewhat
abstractly,
  and it is worth noting that they can be constructed quite
  explicitly.  For example, it follows from what we've done that
  ${\cal Z}(\Sigma_2)$
  is spanned by graphs as shown in Figure \ref{fg:handlebasis},
  where $\lambda$, $\gamma$, and $\delta$ are any representations and
  $f_1$, $f_2$ are any appropriate intertwiners.  Thus ${\cal Z}(\Sigma_2)=
  \bigoplus_{i,j,m=1}^k W_{i,i^*}^j \tnsr W_{m,m^*}^{j^*}$, where
  $W_{i,j}^m$ is the space of intertwiners from $V_{\lambda_i} \tnsr
  V_{\lambda_j}$ to $V_{\lambda_m}$, a space of dimension $N_{i,j}^m$.
\end{itemize}
\end{remark}
\fig
$$\pic{Hndlbd/handlebasis}{50}{-20}{0}{0}$$
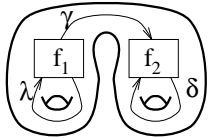
\efig{Spanning labeled ribbon graphs in $H_2$}{fg:handlebasis}

\vspace{30pt}

By necessity the material presented here left off much of interest in
the general field, and gave only a broad sketch of what it did cover.
Below we offer the interested reader some pointers to other
important topics and more details.

Good sources for classical knot theory include the excellent but dated
\cite{Rolfson76} and the more up-to-date \cite{BZ85}.  A  good
elementary survey of the knot polynomials can be found in \cite{LM88}.

Quite a lot is being done with quantum groups, with an eye towards
both  topological and algebraic applications, such as
Kashiwara and Lusztig's canonical bases and the Kazhdan-Lusztig conjectures.
Fortunately, the field has recently benefited from a number of
expository books, including
\cite{Lusztig93,CP94,Kassel94,SS94}, which cover all of this,
and all of the background and details skipped in Sections 2 and 3.
\cite{Baxter82} gives a good account of the integrable models in
statistical mechanics, \cite{Majid90,Fadeev84} connect them to quantum
groups, and  \cite{Kauffman91} presents many ideas of statistical mechanics
from a very knot-theoretic viewpoint.

Much has also been written about the $3$-manifold invariants
constructed in \cite{RT91}.  The $sl_2$ invariant has been constructed
in many different forms and from many different perspectives in
\cite{KM91,Lickorish91,Crane91,BHMV92,BHMV??,Kohno92,Morton92,TW93,Wenzl93,KL94}.
 The square norm of the $3$-manifold invariant associated to any
 modular Hopf algebra can be computed as a certain sum over states on
 a triangulation of the $3$-manifold \cite{TV92,Walker??,Roberts??}.
 There are also a number of $3$-manifold invariants which are
 constructed from the same sorts of data but which appear to be
 different, including \cite{Kuperberg89,Kuperberg??,KR??}.

The unfortunately unpublished \cite{Walker??} gives a thorough and
informative account of the TQFT's we construct, extending the
formalism to cobordisms, with corners (this is a strong version of the
`duality' which Witten uses to get his solution).  It is a good
resource for novices and experts.  Much general imformation can be
found in Quinn's lecture notes on the subject \cite{Quinn95}.  The
other articles in the same volume represent an excellent introduction
to many of the physical and geometric aspects of this subject.

There are other ways to get at the link and $3$-manifold invariants
quite apart from quantum groups.  Affine Lie algebras
\cite{KP84,KW88}, by way of conformal field theory \cite{TK88,Verlinde88},
offer a functor from the tangle category to a sort of intermediate
linear category.  One can still construct the TQFT in this context.
This approach is technically more difficult, because there are
infinite-dimensional algebras and vector spaces to contend with, but is
directly linked to the physics.

Also closer to the physics are efforts to bring the tools
mathematical physics has developed for understanding quantum field
theories to bear on Chern-Simons theory. Good overviews of this
include \cite{Axelrod91,Atiyah90b}.  More recent efforts have focused
on perturbative approaches to Chern-Simons theory, including
\cite{AS92,AS94,Kontsevich??}.  Perturbation theory when applied
to the link invariants is entirely combinatorial, and fits naturally
into the framework of Vassiliev invariants, which were developed out
of purely topological consideration \cite{Vassiliev90,Vassiliev92}.  Good
sources
for the theory of Vassiliev invariants and their relation with the
link invariants in this paper are \cite{BarNatan91,BirmanLin91}, as
well as the excellent expository article which recently appeared in
these pages \cite{Birman93} (it also gives a good introduction to the
link invariants and their place in the history of knot theory).

\newcommand{\etalchar}[1]{$^{#1}$}

\end{document}